\documentstyle{amsppt}
\magnification=\magstep1
\pagewidth{16truecm}
\pageheight{23truecm}
\vcorrection{0truecm}
\hcorrection{0truecm}
\loadbold 
\topmatter
   \title 
         Two permanents \\
         in the universal enveloping algebras
         of the symplectic Lie algebras
   \endtitle
   \author Minoru ITOH \endauthor
   \affil Department of Mathematics and Computer Science, 
          Faculty of Science, \\
          Kagoshima University, Kagoshima 890-0065, Japan \\ 
          itoh\@sci.kagoshima-u.ac.jp \endaffil
   \abstract
   \nofrills{\bf Abstract.}\usualspace
   This paper presents new generators for the center 
   of the universal enveloping algebra 
   of the symplectic Lie algebra.
   These generators are expressed in terms of the column-permanent,
   and it is easy to calculate their eigenvalues 
   on irreducible representations.
   We can regard these generators as the counterpart of
   central elements of  the universal enveloping algebra 
   of the orthogonal Lie algebra given 
   in terms of the column-determinant by A.~Wachi.
   The earliest prototype of all these central elements is 
   the Capelli determinants in the universal enveloping algebra 
   of the general linear Lie algebra.
   \endabstract
\endtopmatter
\leftheadtext{M. Itoh}
\rightheadtext{Two permanents for the symplectic Lie algebras}
%
\document
\define\sgn{\operatorname{sgn}}

\define\diag{\operatorname{diag}}
\define\Pf{\operatorname{Pf}}
\define\Hf{\operatorname{Hf}}
\define\Det{\operatorname{Det}}
\define\per{\operatorname{per}}
\define\Per{\operatorname{Per}}

\define\Mat{\operatorname{Mat}}

\define\Ad{\operatorname{Ad}}

\define\Zgeqzero
{\Bbb{Z}_{{\scriptscriptstyle \geq} \smash{\lower.3ex\hbox{$\scriptstyle 0$}}}}
\define\fischer#1#2{\langle #1 \,|\, #2 \rangle}
\define\iddots{\smash{\lower-.25ex\hbox{.}}\;\!\smash{\lower-.95ex\hbox{.}}\;\!\smash{\lower-1.65ex\hbox{.}}}
\subhead Introduction \endsubhead
In this paper we give new generators for the center 
of the universal enveloping algebra 
of the symplectic Lie algebra $\frak{sp}_N$.
These generators $D_k(u)$ are expressed 
in terms of the ``column-permanent,''
and similar to the Capelli determinants, 
i.e., well-known central elements of 
the universal enveloping algebra
of the general linear Lie algebra $\frak{gl}_N$.
As the key of the Capelli identity,
these Capelli determinants are used 
to analyze the representations of $\frak{gl}_N$
acting via the polarization operators 
(see \cite{Ca1}, \cite{Ca2}, \cite{HU}).
One of the remarkable properties of the Capelli determinants
is that we can  easily calculate their eigenvalues on irreducible representations.
It is also easy to calculate
the eigenvalues of our central elements $D_k(u)$.

On the other hand, it is not so obvious that $D_k(u)$ is actually central 
in the universal enveloping algebra.
This fact can be proved as follows.
In addition to $D_k(u)$,
we consider another central element $D'_k(u)$
expressed in terms of the  ``symmetrized permanent.''
We can easily check that this $D'_k(u)$ is central,
but its eigenvalue is difficult to calculate.
In spite of this difference,
these $D_k(u)$ and $D'_k(u)$ are actually equal.
We will prove this coincidence to show the centrality of $D_k(u)$.
Then, at the same time, we can also see the eigenvalue of $D'_k(u)$.

More directly, our central elements 
are regarded as the counterpart of the central elements 
of the universal enveloping algebra 
of the orthogonal Lie algebra $\frak{o}_N$ 
recently given by A. Wachi \cite{W}
in terms of the ``column-determinant.''
The discussion between $D(u)$ and $D'(u)$ above
can be applied to Wachi's elements (\cite{I4}).

Let us explain the main result precisely.
Let $J \in \Mat_N(\Bbb{C})$ be 
a non-degenerate alternating matrix of size $N$. 
We can realize the symplectic Lie group as the isometry group 
with respect to the bilinear form determined by $J$:
$$
   Sp(J) = \{ g \in GL_N \,|\, {}^t\!g J g = J \}.
$$
The corresponding Lie algebra is expressed as 
$$
   \frak{sp}(J) = \{ Z \in \frak{gl}_N \,|\, {}^t\!Z J + J Z = 0\}.
$$
As generators of this $\frak{sp}(J)$,
we can take $F^{\frak{sp}(J)}_{ij} = E_{ij} - J^{-1} E_{ji} J$,
where $E_{ij}$ is the standard basis of $\frak{gl}_N$.
We introduce the $N \times N$ matrix $F^{\frak{sp}(J)}$
whose $(i,j)$th entry is this generator:
$F^{\frak{sp}(J)} = (F^{\frak{sp}(J)}_{ij})_{1 \leq i,j \leq N}$.
We regard this matrix as an element of $\Mat_N(U(\frak{sp}(J)))$.

In the representation theory,
the case 
$$
   J = J_0 =
   \left(
   \smallmatrix
      & & & & & 1 \\
      & & & & \lower.6ex\hbox{\iddots} & \\
      & & & 1 & & \\
      & & -1 & & & \\
      & \lower.6ex\hbox{\iddots} & & & & \\
      -1 & & & & &
   \endsmallmatrix
   \right)
$$
is important.
Indeed, we can take a triangular decomposition of $\frak{sp}(J_0)$
simply as follows:
$$
   \frak{sp}(J_0) = \frak{n}^- \oplus \frak{h} \oplus \frak{n}^+.
$$
Here $\frak{n}^-$, $\frak{h}$, and $\frak{n}^+$
are the subalgebra of $\frak{sp}(J_0)$ spanned by the elements 
$F^{\frak{sp}(J_0)}_{ij}$ such that
$i>j$, $i=j$, and $i<j$ respectively.
Namely, the entries in the lower triangular part,
in the diagonal part, and in the upper triangular part  
of the matrix $F^{\frak{sp}(J_0)}$
belong to $\frak{n}^-$, $\frak{h}$, and $\frak{n}^+$ respectively.
We call this $\frak{sp}(J_0)$ be the ``split realization''
of the symplectic Lie algebras.

The following is the main theorem of this paper:

\proclaim {Theorem~A} {\bf }
   \sl The following element is central in $U(\frak{sp}(J_0))$
   for any $u \in \Bbb{C}$:
   $$
      D_k(u) = \sum_{1 \leq \alpha_1 \leq \cdots \leq \alpha_k \leq N} 
      \frac{1}{\alpha!} 
      \per(\widetilde{F}^{\frak{sp}(J_0)}_{\alpha} + u \boldkey{1}_{\alpha}
      - \boldkey{1}_{\alpha} 
        \diag(\tfrac{k}{2}-1,\tfrac{k}{2}-2,\ldots,-\tfrac{k}{2})).
   $$
\endproclaim

The notation is as follows.
First, the symbol ``$\per$'' means the ``column-permanent.''
Namely, for an $N \times N$ matrix $Z = (Z_{ij})$, we put
$$
   \per Z = \sum_{\sigma\in\frak{S}_N} 
   Z_{\sigma(1)1} Z_{\sigma(2)2} \cdots Z_{\sigma(N)N},
$$
even if the entries $Z_{ij}$ are non-commutative.
Secondly, $\widetilde{F}^{\frak{sp}(J_0)}$ means the matrix
$$
      \widetilde{F}^{\frak{sp}(J_0)} = 
         F^{\frak{sp}(J_0)} - \diag(0,\ldots,0,1,\ldots,1).
$$
Here the numbers of $0$'s and $1$'s are equal to $N/2$.
Thirdly, $\boldkey{1}$ means the unit matrix.
Moreover, for a matrix $Z =(Z_{ij})$ 
and a non-decreasing sequence $\alpha = (\alpha_1,\ldots,\alpha_k)$,
we denote the matrix $(Z_{\alpha_i \alpha_j})_{1 \leq i,j \leq k}$ 
by $Z_{\alpha}$.
Finally, we put $\alpha! = m_1! \cdots m_N!$,
where $m_1, \ldots, m_N$ are 
the multiplicities of $\alpha = (\alpha_1,\ldots,\alpha_k)$:
$$
    \alpha = (\alpha_1,\ldots,\alpha_k) 
           = (\oversetbrace{m_1}\to{1,\ldots,1},
              \oversetbrace{m_2}\to{2,\ldots,2},\ldots,
              \oversetbrace{m_N}\to{N,\ldots,N}).
$$

This central element $D_k(u)$ is remarkable,
because we can easily calculate its eigenvalue
on irreducible representations of $\frak{sp}(J_0)$ (Theorem~4.4).
However, 
the centrality of $D_k(u)$ (namely Theorem~A) is not so obvious.

To prove Theorem~A, 
we consider another central element of $U(\frak{sp}(J_0))$:
$$
   D'_k(u) = \Per_k (F + u \boldkey{1} \,;\, 
        \tfrac{k}{2}-1, \tfrac{k}{2}-2, \ldots, -\tfrac{k}{2}+1, 0)
$$
Here the symbol ``$\Per_k$'' means the ``symmetrized permanent.''
Namely we put
$$
   \Per_k(Z \,;\, a_1,\ldots,a_k) 
   = \frac{1}{k!} \sum_{1 \leq \alpha_1 \leq \cdots \leq \alpha_k \leq N} 
   \sum_{\sigma,\sigma'\in\frak{S}_N}
   \frac{1}{\alpha!}
   Z_{\alpha_{\sigma(1)}\alpha_{\sigma'(1)}}(a_1) \cdots 
   Z_{\alpha_{\sigma(k)}\alpha_{\sigma'(k)}}(a_k)
$$
with $Z_{ij}(a) = Z_{ij} + \delta_{ij}a$.
From an invariance of this ``$\Per_k$'' (Proposition~1.9),
the centrality of $D'_k(u)$ is almost obvious.
However, its eigenvalue is difficult to calculate.

These $D_k(u)$ and $D'_k(u)$ are actually equal:

\proclaim {Theorem~B} {\bf }
   \sl We have $D_k(u) = D'_k(u)$.
\endproclaim

The centrality of $D_k(u)$ (namely Theorem~A) and 
the eigenvalue of $D'_k(u)$ are both immediate from this Theorem~B.
Namely, proving Theorem~B directly,
we can settle these two problems at the same time.

Our elements $D_k(u)$ can be regarded as the counterpart of 
the central elements of $U(\frak{o}_N)$ recently given by A. Wachi \cite{W}.
Wachi's elements are also expressed in two different ways.
The first expression $C_k(u)$ 
is given in terms of the ``column-determinant,''
and we can easily calculate its eigenvalue under this expression.
On the other hand, 
the second expression $C'_k(u)$ is given 
in terms of the ``symmetrized determinant,''
and the centrality is almost obvious under this expression.
Proving the coincidence $C_k(u) = C'_k(u)$ directly,
we can settle the following two problems:
(i) the centrality of $C_k(u)$, and 
(ii) the calculation of the eigenvalue of $C'_k(u)$.
See Section~3 and \cite{I4} for the details. 

We note some central elements related to these elements.
First, Wachi's element $C_k(u)$ is equal to the central elements
of $U(\frak{o}_N)$ 
given in terms of the Sklyanin determinant in~\cite{M}
(see also \cite{MNO}).
This coincidence is seen by comparing their eigenvalues.
Moreover, $C_k(0)$ and $D_k(0)$ are equal to 
the central elements defined by eigenvalues in~\cite{MN}.
For these elements,
Capelli type identities are given.

The symmetrized permanent was also introduced 
to give Capelli type identities for reductive dual pairs.
See \cite{I2} and \cite{I3} for these Capelli type identities 
in terms of the symmetrized determinant 
and the symmetrized permanent.
These identities are closely related to $C_k(u)$ and $D_k(u)$.

\vskip .10in

The author is grateful to Professors T\^oru Umeda and Akihito Wachi
for the fruitful discussions.

\vskip .20in

\subhead 1. Capelli type elements for the general linear Lie algebras 
\endsubhead
In this section, we recall the Capelli determinant, 
a famous central element of $U(\frak{gl}_N)$ 
essentially given in \cite{Ca1}.
We also recall its generalization in terms of minors given in \cite{Ca2} 
and its analogue in terms of permanents due to M. Nazarov \cite{N}.
These are the prototypes of the main objects of this paper
and Wachi's elements.

\vskip .10in

\noindent
{\bf 1.1.}
First let us recall the Capelli determinant.
Let $E_{ij}$ be the standard basis of $\frak{gl}_N$,
and consider the matrix 
$E = (E_{ij})_{1 \leq i,j \leq N}$ in 
$\Mat_N(U(\frak{gl}_N))$.
The following ``Capelli determinant'' in $U(\frak{gl}_N)$ is well known 
as the key of the Capelli identity
(\cite{Ca1}, \cite{H}, \cite{U1}):
$$
   C^{\frak{gl}_N}(u) = \det(E + u \boldkey{1} + \diag(N-1,N-2,\ldots,0)).
$$
Here the symbol ``$\det$'' means the ``column-determinant.''
Namely, for $N \times N$ matrix $Z = (Z_{ij})$,
we put
$$
   \det Z = \sum_{\sigma \in \frak{S}_N} \sgn(\sigma) 
   Z_{\sigma(1)1} Z_{\sigma(2)2} \cdots Z_{\sigma(N)N}.
$$
Here each $Z_{ij}$ is an element 
of a (non-commutative) associative $\Bbb{C}$-algebra $\Cal{A}$.
This $C^{\frak{gl}_N}_{\det}(u)$ is known to be central:

\proclaim {Theorem~1.1} {\bf }
   \sl The element $C^{\frak{gl}_N}(u)$ is central in $U(\frak{gl}_N)$
   for any $u \in \Bbb{C}$.
\endproclaim

The eigenvalue of this Capelli determinant 
on irreducible representations is easily calculated:

\proclaim {Theorem~1.2} {\bf }
   \sl For the irreducible representation $\pi^{\frak{gl}_N}_{\lambda}$ 
   of $\frak{gl}_N$ determined by the partition 
   $\lambda = (\lambda_1,\ldots,\lambda_N)$,
   we have
   $$
      \pi^{\frak{gl}_N}_{\lambda}(C^{\frak{gl}_N}(u)) 
       = (u + l_1) \cdots (u + l_N).
   $$
   Here we put $l_i = \lambda_i + N - i$.
\endproclaim

This is immediate from the definition of the column-determinant
and the triangular decomposition
$$
   \frak{gl}_N = \frak{n}^- \oplus \frak{h} \oplus \frak{n}^+.
\tag 1.1
$$
Here $\frak{n}^-$, $\frak{h}$, and $\frak{n}^+$ are
the subalgebras of $\frak{gl}_N$ spanned by the elements $E_{ij}$
such that $i > j$, $i = j$, and $i < j$ respectively.
Namely the entries in the lower triangular part, in the diagonal part,
and in the upper triangular part of $E$ belong to 
$\frak{n}^-$, $\frak{h}$, and $\frak{n}^+$ respectively.
Considering the action of $C^{\frak{gl}_N}(u)$ 
to the highest weight vector,
we can easily check Theorem~1.2.

We can rewrite this Capelli determinant 
in terms of the ``symmetrized determinant'' as follows:

\proclaim {Theorem~1.3} {\bf }
   \sl We have
   $$
      \det(E + u \boldkey{1} + \diag(N-1,N-2,\ldots,0)) 
      = \Det(E + u \boldkey{1} \,;\, N-1,N-2,\ldots,0).
   $$
\endproclaim

Here the symbol ``$\Det$'' means the ``symmetrized determinant.''
Namely, for an $N \times N$ matrix $Z = (Z_{ij})$, 
we put
$$
   \Det Z = \frac{1}{N!} \sum_{\sigma, \sigma' \in \frak{S}_N} 
   \sgn(\sigma) \sgn(\sigma') 
   Z_{\sigma(1)\sigma'(1)} Z_{\sigma(2)\sigma'(2)} 
   \cdots Z_{\sigma(N)\sigma'(N)}.
$$
Moreover, for $N$ parameters $a_1,\ldots,a_N \in \Bbb{C}$,
we put
$$\align
   & \Det(Z \,;\, a_1,\ldots,a_N) \\
   & \qquad = \frac{1}{N!} \sum_{\sigma, \sigma' \in \frak{S}_N} 
   \sgn(\sigma) \sgn(\sigma') 
   Z_{\sigma(1)\sigma'(1)}(a_1) Z_{\sigma(2)\sigma'(2)}(a_2) 
   \cdots Z_{\sigma(N)\sigma'(N)}(a_N)
\endalign
$$
with $Z_{ij}(a) = Z_{ij} + \delta_{ij} a$.
It is obvious that $\Det Z$ is equal to the usual determinant,
if the entries are commutative.
This non-commutative determinant ``$\Det$'' is useful to construct 
central elements in $U(\frak{gl}_N)$.
Indeed, we have the following.

\proclaim {Proposition~1.4} {\bf }
   \sl For any $a_1,\ldots,a_N \in \Bbb{C}$,
   the determinant 
   $$
      \Det(E \,;\, a_1,\ldots,a_N)
   $$
   is invariant under the adjoint action of $GL_N(\Bbb{C})$,
   and hence this is central in $U(\frak{gl}_N)$.
\endproclaim

This is immediate from the following two lemmas:

\proclaim {Lemma~1.5} {\bf }
   \sl The symmetrized determinant is invariant 
   under the conjugation by $g \in GL_N(\Bbb{C})$:
   $$
      \Det(g Z g^{-1} \,;\, a_1,\ldots,a_N)
      = \Det(Z \,;\, a_1,\ldots,a_N).
   $$
   Here $Z$ is an arbitrary $N \times N$ matrix whose entries are
   elements of an associative $\Bbb{C}$-algebra $\Cal{A}$.
\endproclaim

\proclaim {Lemma~1.6} {\bf }
   \sl The matrix $E$  satisfies the relation 
   $$
      \Ad(g) E = {}^tg \cdot E \cdot {}^t\!g^{-1}
   $$
   for any $g \in GL_N(\Bbb{C})$.
   Here $\Ad(g) E$ means the matrix $(\Ad(g)E_{ij})_{1 \leq i,j \leq N}$.
\endproclaim

Lemma~1.6 can be checked by a direct calculation.
Lemma~1.5 is also easy from the expression of ``$\Det$''
in the framework of the exterior calculus.
See \cite{I1} for the details (cf. Section~2 of this paper).

For convenience, we consider the symbol $\natural_k = (k-1,k-2,\ldots,0)$.
Then both sides of Theorem~1.3 can be expressed simply as
$$
   C^{\frak{gl}_N}(u) = \det(E + u \boldkey{1} + \diag \natural_N), \qquad
   C^{\prime \frak{gl}_N}(u) = \Det(E + u \boldkey{1} \,;\, \natural_N).
$$
These two expressions play contrast roles.
Indeed, it is not so easy to calculate the eigenvalue 
of $C^{\prime \frak{gl}_N}(u)$ directly,
but the centrality of~$C^{\prime \frak{gl}_N}(u)$ 
is immediate from Proposition~1.4, because 
$\Det(E + u \boldkey{1} \,;\, \natural_N)
= \Det(E \,;\, u 1_N + \natural_N)$.
Here $u 1_N + \natural_N$ means the linear combination of the two  vectors
$1_N = (1,\ldots,1)$ and $\natural_N$ in~$\Bbb{C}^N$.
Namely we put $u 1_N + \natural_N = (u+N-1,u+N-2,\ldots,u)$.

Using Theorem~1.3, we can settle the following two problems
at the same time:
(i) the centrality of $C^{\frak{gl}_N}(u)$, and
(ii) the calculation of the eigenvalue of $C^{\prime \frak{gl}_N}(u)$.
Indeed, as seen above,
the eigenvalue of $C^{\frak{gl}_N}(u)$
and the centrality of $C^{\prime\frak{gl}_N}(u)$ are almost obvious.

The proof of Theorem~1.3 will be given in Section~2.

\vskip .10in

\noindent
{\bf 1.2.}
Next we recall some generalizations of the Capelli determinant.
First, we put
$$
   C^{\frak{gl}_N}_k(u)
   = \sum_{1 \leq \alpha_1 < \cdots < \alpha_k \leq N} 
     \det(E_{\alpha} + u \boldkey{1} + \diag \natural_k).
$$
Here we denote by $Z_{\alpha}$ 
the submatrix $(Z_{\alpha_i \alpha_j})_{1 \leq i,j \leq k}$ 
of the matrix $Z =(Z_{ij})$.
Obviously we have $C^{\frak{gl}_N}_N(u) = C^{\frak{gl}_N}(u)$.
This element $C^{\frak{gl}_N}_k(u)$ is also central in $U(\frak{gl}_N)$ 
for any $u \in \Bbb{C}$,
and known by the name of the ``Capelli elements of degree $k$.''

Moreover we consider the element 
$$
   D^{\frak{gl}_N}_k(u)
   = \sum_{1 \leq \alpha_1 \leq \cdots \leq \alpha_k \leq N}
   \frac{1}{\alpha!} 
   \per(E_{\alpha} + u \boldkey{1}_{\alpha} 
        - \boldkey{1}_{\alpha} \diag \natural_k)
$$
due to Nazarov \cite{N}.
Here the symbol ``$\per$'' means the ``column-permanent.''
Namely, for any $N \times N$ matrix $Z = (Z_{ij})$,
we put
$$
   \per Z = \sum_{\sigma\in\frak{S}_N} Z_{\sigma(1)1} \cdots Z_{\sigma(N)N}.
$$
Moreover, we put $\alpha! = m_1! \cdots m_N!$,
where $m_1, \ldots, m_N$ are 
the multiplicities of $\alpha = (\alpha_1,\ldots,\alpha_k)$:
$$
    \alpha = (\alpha_1,\ldots,\alpha_k) 
           = (\oversetbrace{m_1}\to{1,\ldots,1},
              \oversetbrace{m_2}\to{2,\ldots,2},\ldots,
              \oversetbrace{m_N}\to{N,\ldots,N}).
$$
Note that $Z_{\alpha} = (Z_{\alpha_i \alpha_j})_{1 \leq i,j \leq k}$
is not a submatrix of $Z$ in general,
because $\alpha$ has some multiplicities.
This $D^{\frak{gl}_N}_k(u)$ 
is also central in $U(\frak{gl}_N)$ for any $u \in \Bbb{C}$.

We can easily calculate the eigenvalues of these elements
$C^{\frak{gl}_N}_k(u)$ and $D^{\frak{gl}_N}_k(u)$.
The proof is almost the same as that of Theorem~1.2:

\proclaim {Proposition~1.7} {\bf }
   \sl For the irreducible representation $\pi = \pi^{\frak{gl}_N}_{\lambda}$
   of $\frak{gl}_N$ determined by the partition 
   $\lambda = (\lambda_1,\ldots,\lambda_N)$, we have
   $$\align
      \pi(C^{\frak{gl}_N}_k(u)) 
      & = \sum_{1 \leq \alpha_1 < \cdots < \alpha_k \leq N} 
      (u+\lambda_{\alpha_1}+k-1) (u+\lambda_{\alpha_2}+k-2) 
      \cdots (u+\lambda_{\alpha_k}), \\
      \pi(D^{\frak{gl}_N}_k(u)) 
      & = \sum_{1 \leq \alpha_1 \leq \cdots \leq \alpha_k \leq N} 
      (u+\lambda_{\alpha_1}-k+1) (u+\lambda_{\alpha_2}-k+2) 
      \cdots (u+\lambda_{\alpha_k}).
   \endalign
   $$
\endproclaim

We can rewrite these elements in terms of the ``symmetrized determinant''
and the ``symmetrized permanent'':

\proclaim {Theorem~1.8} {\bf }
   \sl We have 
   $$\align
     \sum_{1 \leq \alpha_1 < \cdots < \alpha_k \leq N} 
     \det(E_{\alpha} + u \boldkey{1} + \diag \natural_k)
      & = \Det_k(E + u \boldkey{1} \,;\, \natural_k), \\
     \sum_{1 \leq \alpha_1 \leq \cdots \leq \alpha_k \leq N}
     \frac{1}{\alpha!} 
     \per(E_{\alpha} + u \boldkey{1}_{\alpha} 
        - \boldkey{1}_{\alpha} \diag \natural_k)
      & = \Per_k(E + u \boldkey{1} \,;\, -\natural_k).
   \endalign
   $$
\endproclaim

Here $\Det_k$ and $\Per_k$ are defined as follows.
First we put
$$
   \Per Z = \frac{1}{N!} \sum_{\sigma,\sigma'\in\frak{S}_N} 
     Z_{\sigma(1)\sigma'(1)} \cdots Z_{\sigma(k)\sigma'(k)}.
$$
Noting this, we put
$$\align
   \Det_k(Z) 
   & = \sum_{1 \leq \alpha_1 < \cdots < \alpha_k \leq N}
       \Det Z_{\alpha} \\
   & = \sum_{1 \leq \alpha_1 < \cdots < \alpha_k \leq N}
     \frac{1}{k!}
     \sum_{\sigma,\sigma'\in\frak{S}_k}
     \sgn(\sigma)\sgn(\sigma')
     Z_{\alpha_{\sigma(1)} \alpha_{\sigma'(1)}}\cdots 
     Z_{\alpha_{\sigma(k)} \alpha_{\sigma'(k)}}, \\
   \Per_k(Z) 
   & = \sum_{1 \leq \alpha_1 \leq \cdots \leq \alpha_k \leq N}
       \frac{1}{\alpha!} \Per Z_{\alpha} \\
   & = \sum_{1 \leq \alpha_1 \leq \cdots \leq \alpha_k \leq N}
     \frac{1}{\alpha!}\frac{1}{k!} \sum_{\sigma,\sigma'\in\frak{S}_k}
     Z_{\alpha_{\sigma(1)} \alpha_{\sigma'(1)}}\cdots 
     Z_{\alpha_{\sigma(k)} \alpha_{\sigma'(k)}}.
\endalign
$$
Moreover, for $k$ parameters $a_1,\ldots,a_k \in \Bbb{C}$, we put
$$\align
   & \Det_k(Z \,;\, a_1,\ldots,a_k) \\
   & \qquad = \sum_{1 \leq \alpha_1 < \cdots < \alpha_k \leq N} 
     \frac{1}{k!} \sum_{\sigma,\sigma'\in\frak{S}_k} 
     \sgn(\sigma)\sgn(\sigma')
     Z_{\alpha_{\sigma(1)} \alpha_{\sigma'(1)}}(a_1) \cdots 
     Z_{\alpha_{\sigma(k)} \alpha_{\sigma'(k)}}(a_k), \\
   & \Per_k(Z \,;\, a_1,\ldots,a_k) \\
   & \qquad = \sum_{1 \leq \alpha_1 \leq \cdots \leq \alpha_k \leq N} 
     \frac{1}{\alpha!} \frac{1}{k!} \sum_{\sigma,\sigma'\in\frak{S}_k} 
     Z_{\alpha_{\sigma(1)} \alpha_{\sigma'(1)}}(a_1) \cdots 
     Z_{\alpha_{\sigma(1)} \alpha_{\sigma'(1)}}(a_k).
\endalign
$$
These $\Det_k$ and $\Per_k$ are invariant under the conjugations:

\proclaim {Proposition 1.9} {\bf }
   \sl For $g \in GL_N(\Bbb{C})$, we have
   $$\align
      \Det_k(g Z g^{-1} \,;\, a_1,\ldots,a_k)
      & = \Det_k(Z \,;\, a_1,\ldots,a_k), \\
      \Per_k(g Z g^{-1} \,;\, a_1,\ldots,a_k)
      & = \Per_k(Z \,;\, a_1,\ldots,a_k).
   \endalign
   $$
\endproclaim

Hence, combining this with Lemma~1.6,
we have the following.

\proclaim {Proposition 1.10} {\bf }
   \sl For arbitrary $a_1,\ldots,a_k \in \Bbb{C}$,
   the two elements
   $$
      \Det_k(E \,;\, a_1,\ldots,a_k), \qquad
      \Per_k(E \,;\, a_1,\ldots,a_k)
   $$
   are invariant under the adjoint action of $GL_N(\Bbb{C})$,
   and hence central in $U(\frak{gl}_N)$.
\endproclaim

Let us denote 
by $C^{\prime \frak{gl}_N}_k(u)$ and $D^{\prime \frak{gl}_N}_k(u)$ 
the right hand sides of Theorem~1.8:
$$\align
   C^{\prime \frak{gl}_N}_k(u)
   & = \Det_k(E + u \boldkey{1} \,;\, \natural_k)
   = \Det_k(E \,;\, u 1_N + \natural_k), \\
   D^{\prime \frak{gl}_N}_k(u)
   & = \Per_k(E + u \boldkey{1} \,;\, -\natural_k)
   = \Per_k(E \,;\, u1_N-\natural_k).
\endalign
$$
These are obviously central  
in $U(\frak{gl}_N)$ for any $u \in \Bbb{C}$.
However it is not so easy to calculate their eigenvalues directly.

Theorem~1.8 settles the following two problems at the same time:
(i) the centralities of $C^{\frak{gl}_N}_k(u)$ and $D^{\frak{gl}_N}_k(u)$,
and
(ii) the calculation of the eigenvalues of  
$C^{\prime \frak{gl}_N}_k(u)$ and~$D^{\prime \frak{gl}_N}_k(u)$.

\vskip .20in

\subhead 2. The proof in the case of the general linear Lie algebras 
\endsubhead
In this section, 
we recall the proofs of Theorems~1.3 and 1.8 given 
in \cite{I1} and \cite{I3}.
In these proofs,
we express determinants and permanents 
in the framework of the exterior algebra and the symmetric tensor algebra.
These calculations are the prototypes of the proof of the main theorem.

\vskip .10in

\noindent
{\bf 2.1.}
First, we recall the proof of Theorem~1.3 given in \cite{I1}.

We can express the column-determinant
in the framework of the exterior algebra as follows.
Let $e_1,\ldots,e_N$ be $N$ anti-commuting formal variables,
which generate the exterior algebra $\Lambda_N = \Lambda(\Bbb{C}^N)$.
Put $\eta_j(u) = \sum_{i=1}^N e_i E_{ij}(u)$
as an element in the extended algebra
 $\Lambda_N \otimes U(\frak{gl}_N)$ 
in which the two subalgebras $\Lambda_N$ and $U(\frak{gl}_N)$ 
commute with each other.
Then, by a direct calculation, we have the following equality
in $\Lambda_N \otimes U(\frak{gl}_N)$:
$$
   \eta_1(a_1) \eta_2(a_2) \cdots \eta_N(a_N)
   = e_1 e_2 \cdots e_N
     \det(E + \diag(a_1,a_2,\ldots,a_N)).
\tag 2.1
$$

The symmetrized determinant is expressed similarly 
by doubling the anti-commuting variables.
Let $e_1,\ldots,e_N, e^*_1,\ldots,e^*_N$ be 
$2N$ anti-commuting formal variables,
which generate the exterior algebra 
$\Lambda_{2N} = \Lambda(\Bbb{C}^N \oplus \Bbb{C}^N)$.
We put $\varXi(u) = \sum_{i,j=1}^N e_i e^*_j E_{ij}(u)$
in $\Lambda_{2N} \otimes U(\frak{gl}_N)$ . 
Then, by a direct calculation, we have
$$
   \varXi(a_1) \varXi(a_2) \cdots \varXi(a_N) 
   = N! e_1 e^*_1 \cdots e_N e^*_N \Det(E \,;\, a_1,a_2,\ldots,a_N).
\tag 2.2
$$

Now we can prove Theorem~1.3 using the commutation relation
$$
   \eta_{i} (a+1) \eta_{j} (a) + \eta_{j} (a+1) \eta_{i} (a) = 0.
\tag 2.3
$$
This relation itself is easy from the relation 
$[E_{ij},E_{kl}] = \delta_{kj} E_{il} - \delta_{il} E_{kj}$.

\vskip .10in

\noindent
{\it Proof of Theorem~{\sl 1.3.}}
Since $\varXi(u) = \sum_{i=1}^N \eta_i (u) e^*_i$, we have
$$\align
   & \varXi(u+N-1) \varXi(u+N-2) \cdots \varXi(u) \\
   & \quad = \sum_{1 \leq i_1,\ldots,i_N \leq N} 
     \eta_{i_1} (u+N-1) e^*_{i_1} 
     \eta_{i_2} (u+N-2) e^*_{i_2} 
     \cdots 
     \eta_{i_N} (u) e^*_{i_N} \\
   & \quad = 
     (-)^{\frac{N(N-1)}{2}}
     \sum_{1 \leq i_1,\ldots,i_N \leq N} 
     \eta_{i_1} (u+N-1) \eta_{i_2} (u+N-2) \cdots \eta_{i_N} (u)
     e^*_{i_1}  e^*_{i_2} \cdots e^*_{i_N}.
\endalign
$$
Here the indices $i_1,\ldots,i_N$  can be regarded as 
a permutation of $1,\ldots,N$, because $e^*_i$'s are anti-commuting.
Moreover the factors~$\eta_{i} (a)$ can be reordered by
using the commutation relation (2.3).
Thus we have
$$\align
   & \varXi(u+N-1) \varXi(u+N-2) \cdots \varXi(u) \\
   & \quad 
   = (-)^{\frac{N(N-1)}{2}} \sum_{\sigma \in \frak{S}_N}
     \eta_{\sigma(1)} (u+N-1)
     \eta_{\sigma(2)} (u+N-2)
     \cdots 
     \eta_{\sigma(N)} (u) 
     e^*_{\sigma(1)} e^*_{\sigma(2)} \cdots e^*_{\sigma(N)} \\
   & \quad 
   = (-)^{\frac{N(N-1)}{2}} N!
     \eta_1 (u+N-1) \eta_2 (u+N-2) \cdots \eta_N (u) 
     e^*_1 e^*_2 \cdots e^*_N.
\endalign
$$
Comparing this equality with (2.1) and (2.2),
we reach to the assertion.
\hfill\qed

\vskip .10in

\noindent
{\it Remark.}
From (2.2), we can see that $\Det(E \,;\, a_1,\ldots,a_N)$ does not 
depend on the order of the parameters $a_1,\ldots,a_N$,
because $\varXi(a_1),\ldots,\varXi(a_N)$ commute with each other.
Indeed $\varXi(u)$ can be expressed as $\varXi(u) = \varXi(0) + u \tau$
with $\tau = \sum_{i=1}^N e_i e^*_i$, 
and this $\tau$ is central in $\Lambda_{2N} \otimes U(\frak{gl}_N)$.

\vskip .10in

\noindent
{\bf 2.2.}
Next we go to the proof of Theorem~1.8.
We only prove the second relation here,
because the proof of the first one is almost the same.

We start with the expressions of our permanents
in the framework of the symmetric tensor algebra.
Let $Z = (Z_{ij})$ be an $N \times N$ matrix
whose entries are elements 
of a (non-commutative) associative $\Bbb{C}$-algebra $\Cal{A}$.
Let $e_1,\ldots,e_N$ be $N$ commutative formal variables,
which generate the symmetric tensor algebra $S_N = S(\Bbb{C}^N)$.
We put $\eta_j = \sum_{i=1}^N e_i Z_{ij}$
as an element in the extended algebra $S_N \otimes \Cal{A}$ 
in which the two subalgebras $S_N$ and $\Cal{A}$ 
commute with each other.
Then, by a direct calculation, we have the relation
$$
   \eta_{\beta}
   = \sum_{1 \leq \alpha_1 \leq \cdots \leq \alpha_k \leq N}
     \frac{1}{\alpha!} e_{\alpha}
     \per(Z_{\alpha\beta})
\tag 2.4
$$
for $1 \leq \beta_1 \leq \cdots \leq \beta_k \leq N$.
Here  $\eta_{\beta}$ and $e_{\alpha}$ denote 
$\eta_{\beta_1} \cdots \eta_{\beta_k}$ 
and $e_{\alpha_1} \cdots e_{\alpha_k}$ respectively,
and $Z_{\alpha\beta}$ means the matrix
$Z_{\alpha\beta} = (Z_{\alpha_i \beta_j})_{1 \leq i,j \leq k}$.
Moreover, by putting $\eta_j(u) = \sum_{i=1}^N e_i Z_{ij}(u)$,
this relation is generalized to
$$
   \eta_{\beta_1}(a_1) \cdots \eta_{\beta_k}(a_k) 
   = \sum_{1 \leq \alpha_1 \leq \cdots \leq \alpha_k \leq N}
     \frac{1}{\alpha!} e_{\alpha} 
      \per(Z_{\alpha\beta} 
           + \boldkey{1}_{\alpha\beta} \diag(a_1,\ldots,a_k)).
\tag 2.5
$$

The symmetrized permanent 
is similarly expressed by doubling the commutative variables.
Let $e_1,\ldots,e_N, e^*_1,\ldots,e^*_N$ be 
$2N$ commutative formal variables,
which generate the symmetric tensor algebra 
$S_{2N} = S(\Bbb{C}^N \oplus \Bbb{C}^N)$.
We put $\varXi = \sum_{i,j=1}^N e_i e^*_j Z_{ij}$ 
in $S_{2N} \otimes \Cal{A}$ . 
Then, we have the relation
$$
   \varXi^{(k)} = 
   \sum \Sb 1 \leq \alpha_1 \leq \cdots \leq \alpha_k \leq N \\
            1 \leq \beta_1 \leq \cdots \leq \beta_k \leq N \endSb
   \frac{1}{\alpha!\beta!} 
   e_{\alpha} e^*_{\beta} \Per(Z_{\alpha\beta}).
\tag  2.6
$$
Here $x^{(k)}$ means the divided power: $x^{(k)} = \frac{1}{k!} x^k$.

It is convenient to consider the bilinear form 
$\langle \cdot \,|\, \cdot \rangle$ on $S_{2N}$
defined by the formula
$$
   \big< e_{\alpha} e^*_{\beta} \,\big|\, e_{\alpha'} e^*_{\beta'} \big>
   = \delta_{\alpha\alpha'} \delta_{\beta\beta'} \alpha! \beta!.
$$
This is known as the ``Fischer inner product.''
Using this, we can rewrite (2.4) and (2.5) to get
$$ 
   \per(Z_{\alpha\beta}) 
   = \big< \eta_{\beta} \,\big|\, e_{\alpha} \big>, \quad 
   \per(Z_{\alpha\beta} 
        + \boldkey{1}_{\alpha\beta} \diag(a_1,\ldots,a_k)) 
   = \big< \eta_{\beta_1}(a_1)\cdots\eta_{\beta_k}(a_k) \,\big|\, 
             e_{\alpha} \big>.
$$
Similarly (2.6) can be rewritten to get
$$
   \Per(Z_{\alpha\beta}) = \big< \varXi^{(k)} \,\big|\, 
   e_{\alpha} e^*_{\beta} \big>.
$$
Moreover, putting $\varXi(u) = \sum_{i,j=1}^N e_i e^*_j Z_{ij}(u)$ 
and $\tau = \sum_{i=1}^N e_i e^*_i$,
we can express $\Per_k$ as 
$$
   \Per_k(Z) = \big< \varXi^{(k)} \,\big|\, \tau^{(k)} \big>, \qquad
   \Per_k(Z \,;\, a_1,\ldots,a_k) 
   = \big< \frac{1}{k!}\varXi(a_1)\cdots\varXi(a_k) \,\big|\, 
             \tau^{(k)} \big>.
$$
These are immediate by noting the relation
$$
   \tau^{(k)} 
   = \sum_{1 \leq \alpha_1 \leq \cdots \leq \alpha_k \leq N}
     \frac{1}{\alpha!} e_{\alpha} e^*_{\alpha}.
$$
We have a similar expression for the column-permanent:
$$
   \per(Z_{\alpha} + u \boldkey{1}_{\alpha} 
   + \boldkey{1}_{\alpha} \diag(a_1,\ldots,a_k)) 
   = 
   \big<
      \eta^{\dagger}_{\alpha_1}(u+a_1) \cdots \eta^{\dagger}_{\alpha_k}(u+a_k)
   \,\big|\, \tau^{(k)} \big>.
$$
Here we put $\eta^{\dagger}_j(u) = \eta_j(u) e^*_j$.

Let us write these expressions simply as
$$
      \Per_k(Z) = \big< \varXi^{(k)} \big>, \qquad
   \Per_k(Z \,;\, a_1,\ldots,a_k) 
   = \big< \frac{1}{k!}\varXi(a_1)\cdots\varXi(a_k) \big>,
\tag 2.7
$$
$$
   \per(Z_{\alpha} + u \boldkey{1}_{\alpha} 
   + \boldkey{1}_{\alpha} \diag(a_1,\ldots,a_k)) 
   = 
   \big<
      \eta^{\dagger}_{\alpha_1}(u+a_1) \cdots \eta^{\dagger}_{\alpha_k}(u+a_k)
   \big>.
\tag 2.8
$$
Here we put
$$
   \big< \varphi \big> 
   = \sum_{k=0}^{\infty} \big< \varphi \,\big|\, \tau^{(k)}\big>
$$
for $\varphi \in S_{2N}$.
Note here that the sum is actually finite.

\vskip .10in

\noindent
{\it Remark.}
We can see that $\Per_k(Z \,;\, a_1,\ldots,a_k)$
does not depend on the order the parameters $a_1,\ldots,a_k$,
because $\varXi(a_1),\ldots,\varXi(a_N)$ commute with each other.
Indeed, we can express $\varXi(u)$ as $\varXi(0) + u \tau$, 
and $\tau$ is central in $S_{2N} \otimes \Cal{A}$.
Similarly, $\Det_k(Z \,;\, a_1,\ldots,a_k)$
does not depend on the order the parameters.

\vskip .10in

Using these expressions,
we can prove the second relation of Theorem~1.8 as follows:

\vskip .10in

\noindent
{\it Proof of the second relation of Theorem~{\sl 1.8.}}
We put $\eta_i (u) = \sum_{j=1}^N e_j E_{ij}(u)$.
The commutation relation 
$$
   \eta_{i} (a) \eta_{j} (a+1) - \eta_{j} (a) \eta_{i} (a+1) = 0
$$
is easy from the relation 
$[E_{ij},E_{kl}] = \delta_{kj} E_{il} - \delta_{il} E_{kj}$.
In particular, $\eta^{\dagger}_i (u) = \eta_i (u) e^*_i$
satisfies the relation
$$
   \eta^{\dagger}_{i} (a) \eta^{\dagger}_{j} (a+1) 
   - \eta^{\dagger}_{j} (a) \eta^{\dagger}_{i} (a+1) = 0.
\tag 2.9
$$
Since $\varXi(u) = \sum_{i,j=1}^N e_i e^*_j E_{ij}(u)$ 
is written as $\varXi(u) = \sum_{i=1}^N \tilde{\eta}_i (u)$,
we have
$$\align
   & \varXi(u-k+1) \varXi(u-k+2) 
     \cdots \varXi(u) \\
   & \qquad = \sum_{1 \leq i_1,\ldots,i_k \leq N} 
     \eta^{\dagger}_{i_1} (u-k+1) 
     \eta^{\dagger}_{i_2} (u-k+2)
     \cdots 
     \eta^{\dagger}_{i_k} (u).
\endalign
$$
The factors~$\eta^{\dagger}_{i} (a)$ can be reordered by
using the commutation relation (2.9).
Thus we have
$$\align
   & \varXi(u-k+1) \varXi(u-k+2) 
     \cdots \varXi(u) \\
   & \qquad 
   = \sum_{1 \leq \alpha_1 \leq \cdots \leq \alpha_k \leq N} 
     \frac{k!}{\alpha!}
     \eta^{\dagger}_{\alpha_1} (u-k+1) 
     \eta^{\dagger}_{\alpha_2} (u-k+2) 
     \cdots 
     \eta^{\dagger}_{\alpha_k} (u).
\endalign
$$
Comparing this equality with (2.7) and (2.8), we reach to the assertion.
\hfill\qed

\vskip .10in

Next, let us prove Proposition~1.9.
This is an application of the following lemma,
an elementary fact for the ``Fischer inner product''
$\left< \cdot \,|\, \cdot \right>$:

\proclaim {Lemma~2.1} {\bf }
   \sl Consider the standard action of $g \in GL_{2N}$ 
   on the vector space $\Bbb{C}^N \oplus \Bbb{C}^N$,
   which is naturally extended to an automorphism 
   of $S_{2N}$.
   In this situation, we have
   $$
      \big< \varphi \,\big|\, \varphi' \big>
      =  \big< g (\varphi) \,\big|\, {}^t\!g^{-1}(\varphi') \big>
   $$
   for $\varphi, \varphi' \in S_{2N}$.
\endproclaim

\noindent
{\it Proof of the second relation of Proposition~{\sl 1.9.}}
For $g \in GL_N$, we put
$$
   h_g = \diag(g, {}^t\!g^{-1}) 
   = \pmatrix
     g & 0 \\
     0 & {}^t\!g^{-1}
     \endpmatrix
   \in GL_{2N},
$$
and consider its natural action on $S_{2N} \otimes \Cal{A}$.
By direct calculations, we have the relations
$h_g(\varXi_{Z}(u)) = \varXi_{gZg^{-1}}(u)$ 
and ${}^t\!h_g^{-1}(\tau) = \tau$
for $\varXi_Z(u) = \sum_{i,j = 1}^N e_i e^*_j Z_{ij}(u)$ and 
$\tau = \sum_{i=1}^N e_i e^*_i$.
Since $h_g$ and ${}^t\!h_g^{-1}$ are 
automorphisms of $S_{2N} \otimes \Cal{A}$,
we have
$$\align
   \big< \varXi_{gZg^{-1}}(u_1) \cdots \varXi_{gZg^{-1}}(u_k) \,\big|\,
           \tau^{(k)} \big>
   & = \big< h_g (\varXi_Z(u_1) \cdots h_g \varXi_Z(u_k)) \,\big|\,
               {}^t\!h_g^{-1}(\tau)^{(k)} \big> \\
   & = \big< \varXi_Z(u_1) \cdots \varXi_Z(u_k) \,\big|\,
               \tau^{(k)} \big>.
\endalign
$$
Here, we used Lemma~2.1 for the second equality.
By (2.7) this implies our assertion.
\hfill\qed

\vskip .10in

The relations for determinants in Theorem~1.8 and Proposition~1.9
can be proved similarly
by considering the exterior algebra 
instead of the symmetric tensor algebra
(see \cite{I3}).

\vskip .20in

\subhead 3. The case of the orthogonal Lie algebras 
\endsubhead
Before going to the main result 
in the case of the symplectic Lie algebra $\frak{sp}_N$,
we recall the case of the orthogonal Lie algebra $\frak{o}_N$.
In this case, two analogues of the Capelli determinant are known.
One was given by R. Howe and T. Umeda \cite{HU},
and the other was recently given by A. Wachi~\cite{W}.

\vskip .10in


\noindent
{\bf 3.1.}
First we see the general realization of $\frak{o}_N$.
Let $S \in \Mat_N(\Bbb{C})$ be a nondegenerate symmetric matrix of size $N$.
We can realize the orthogonal Lie group
as the isometry group with respect to the bilinear form
determined by $S$:
$$
   O(S) = \{ g \in GL_N \,|\, {}^t\!g S g = S \}.
$$
The corresponding Lie algebra is expressed as
$$
   \frak{o}(S) = \{ Z \in \frak{gl}_N \,|\, {}^t\!Z S + S Z = 0 \}.
$$
As generators of this $\frak{o}(S)$,
we can take $F^{\frak{o}(S)}_{ij} = E_{ij} - S^{-1} E_{ji} S$,
where $E_{ij}$ is the standard basis of $\frak{gl}_N$.
We consider the $N \times N$ matrix 
$F^{\frak{o}(S)} = (F^{\frak{o}(S)}_{ij})_{1 \leq i,j \leq N}$ 
whose $(i,j)$th entry is this generator $F^{\frak{o}(S)}_{ij}$.
By a direct calculation, 
this $F^{\frak{o}(S)}$ satisfies the following relation:

\proclaim {Lemma~3.1} {\bf }
   \sl For any $g \in O(S)$, we have 
   $$
      \Ad(g) F^{\frak{o}(S)} 
      = {}^t\!g \cdot F^{\frak{o}(S)} \cdot {}^t\!g^{-1}.
   $$
   Here $\Ad(g) F^{\frak{o}(S)}$ means the matrix 
   $(\Ad(g) F^{\frak{o}(S)}_{ij})_{1 \leq i,j \leq N}$.
\endproclaim

Combining this with Proposition~1.9, we have the following.

\proclaim {Proposition~3.2} {\bf }
   \sl The two elements
   $$
      \Det_k(F^{\frak{o}(S)} \,;\, a_1,\ldots,a_k), \qquad
      \Per_k(F^{\frak{o}(S)} \,;\, a_1,\ldots,a_k)
   $$ 
   are invariant 
   under the adjoint action of $O(S)$,
   and in particular these are central in~$U(\frak{o}(S))$:
\endproclaim

Thus,
as in the case of $\frak{gl}_N$,
the symmetrized determinant and the symmetrized permanent 
are useful to construct central elements of $U(\frak{o}(S))$.
On the other hand, unfortunately,
it seems not 
easy to construct central elements of $U(\frak{o}(S))$
similarly using the column-determinant or the column-permanent
at least for general $S$.

However, for some special $S$ ($S=\boldkey{1}$
and $S = S_0 = (\delta_{i,N+1-j})_{1 \leq i,j \leq N}$),
we have analogues of the Capelli determinant
expressed in terms of the column-determinant.

\vskip .10in


\noindent
{\bf 3.2.}
First, let us consider the case that 
$S$ is equal to the unit matrix~$\boldkey{1}$.
Namely we consider the Lie algebra consisting of all alternating matrices:
$$
   \frak{o}(\boldkey{1}) = \{ Z \in \frak{gl}_N \,|\, Z + {}^t\!Z = 0 \}.
$$
In this case, R. Howe and T. Umeda gave an analogue of the Capelli determinant
in terms of the column-determinant:

\proclaim {Theorem~3.3 (\cite{HU})} {\bf }
   \sl The following element is central in $U(\frak{o}(\boldkey{1}))$ 
   for any $u \in \Bbb{C}$:
   $$
      C^{\frak{o}(\boldkey{1})}(u) 
      = \det(F^{\frak{o}(\boldkey{1})} + u \boldkey{1} + \diag\natural_N).
   $$
\endproclaim

As in the case of  $\frak{gl}_N$,
we can rewrite this in terms of the symmetrized determinant:

\proclaim {Theorem~3.4 (\cite{IU})} {\bf }
   \sl We have
   $$
      \det(F^{\frak{o}(\boldkey{1})} 
           + u \boldkey{1} + \diag\natural_N)
      = \Det(F^{\frak{o}(\boldkey{1}_N)} + u \boldkey{1} \,;\, \natural_N).
   $$
\endproclaim

Theorem~3.3 is immediate from this Theorem~3.4.
Indeed, by Proposition~3.2,
$$
   C^{\prime\frak{o}(\boldkey{1})}(u) 
   = \Det(F^{\frak{o}(\boldkey{1})} + u \boldkey{1} \,;\, \natural_N)
   = \Det(F^{\frak{o}(\boldkey{1})} \,;\, u 1_N + \natural_N)
$$
is central in $U(\frak{o}(\boldkey{1}))$ for any $u \in \Bbb{C}$.

\vskip .10in

\noindent
{\it Remarks.}
(1) As in the case of $\frak{gl}_N$, 
we have the following generalization of $C^{\frak{o}(\boldkey{1})}(u)$:
$$
    C^{\frak{o}(\boldkey{1})}_k(u) 
    = \sum_{1 \leq \alpha_1 < \cdots < \alpha_k \leq N}
    \det(F^{\frak{o}(\boldkey{1})}_{\alpha} + u \boldkey{1} 
         + \diag\natural_k).
$$
This can be rewritten in terms of ``$\Det_k$'' as 
$$
   \sum_{1 \leq \alpha_1 < \cdots < \alpha_k \leq N}
   \det(F^{\frak{o}(\boldkey{1})}_{\alpha} 
        + u \boldkey{1} + \diag\natural_k)
   = \Det_k(F^{\frak{o}(\boldkey{1}_N)} + u \boldkey{1} \,;\, \natural_k).
$$

\vskip .10in

\noindent
(2)
These elements are quite similar to 
the Capelli elements $C^{\frak{gl}_N}_k(u)$.
However, it is not easy to calculate their eigenvalues.
Indeed, for this realization $\frak{o}(\boldkey{1})$,
we can not take its triangular decomposition so simply as~(1.1).

\vskip .10in


\noindent
{\bf 3.3.}
Next, we see the case $S = S_0 = (\delta_{i,N+1-j})$,
namely we consider the split realization of the orthogonal Lie algebra:
$$
   \frak{o}(S_0) 
   = \{ Z = (Z_{ij}) \in \frak{gl}_N \,|\, Z_{ij} + Z_{N+1-j, N+1-i} = 0 \}.
$$
A central element of $U(\frak{o}(S_0))$ was recently given
in terms of the column-determinant:

\proclaim {Theorem~3.5 (\cite{W})} {\bf }
   \sl The element
   $$
      C^{\frak{o}(S_0)}(u)
      = \det(F^{\frak{o}(S_0)} + u \boldkey{1} 
        + \diag \tilde{\natural}_N)
   $$
   is central in $U(\frak{o}(S_0))$ for any $u \in \Bbb{C}$.
   Here $\tilde{\natural}_N$ is the following sequence of length $N$:
   $$
      \tilde{\natural}_N =
      \cases
         (\tfrac{N}{2}-1, \tfrac{N}{2}-2, \ldots, 0, 
             0, \ldots, -\tfrac{N}{2}+1), & \text{$N$: even}, \\
         (\tfrac{N}{2}-1, \tfrac{N}{2}-2, \ldots, \tfrac{1}{2}, 0, 
             -\tfrac{1}{2}, \ldots, -\tfrac{N}{2}+1), & \text{$N$: odd}.
      \endcases
   $$
\endproclaim

The proof of this theorem is not so easy
(Wachi showed the commutativity with the generators of $\frak{o}(S_0)$
by employing the exterior calculus).
On the other hand, we can easily calculate its eigenvalue:

\proclaim {Theorem~3.6 (\cite{W})} {\bf }
   \sl Let $\pi^{\frak{o}(S_0)}_{\lambda}$ be the irreducible representation 
   of $\frak{o}(S_0)$ determined by the partition 
   $\lambda = (\lambda_1,\ldots,\lambda_{[N/2]})$,
   where $[N/2]$ means the greatest integer not exceeding $N/2$.
   Then we have
   $$
      \pi^{\frak{o}(S_0)}_{\lambda}(C^{\frak{o}(S_0)}(u))
      = \cases
        (u^2 - l_1^2)(u^2 - l_2^2) \cdots (u^2 - l_{N/2}^2), 
        & \qquad \text{$N$: even}, \\
        u (u^2 - l_1^2)(u^2 - l_2^2) \cdots (u^2 - l_{[N/2]}^2), 
        & \qquad \text{$N$: odd}.
        \endcases
   $$
   Here we put $l_i = \lambda_i + N/2 -i$.
\endproclaim

The proof is almost the same as that of Theorem~1.2.
Namely this is easy from the definition of the column-determinant
and the triangular decomposition
$$
   \frak{o}(S_0) = \frak{n}^- \oplus \frak{h} \oplus \frak{n}^+.
$$
Here $\frak{n}^-$, $\frak{h}$, and $\frak{n}^+$ are 
the subalgebras of $\frak{o}(S_0)$
spanned by the elements $F^{\frak{o}(S_0)}_{ij}$ such that
$i>j$, $i=j$, and $i<j$ respectively.
Namely, the entries in the lower triangular part,
in the diagonal part, and in the upper triangular part
of the matrix $F^{\frak{o}(S_0)}$ belong to 
$\frak{n}^-$, $\frak{h}$, and $\frak{n}^+$ respectively.

We can also rewrite $C^{\frak{o}(S_0)}(u)$ 
in terms of the symmetrized determinant.
To see this we put 
$$
   C^{\prime \frak{o}(S_0)}(u)
   = \Det(F^{\frak{o}(S_0)} + u \boldkey{1} \,;\, \tilde{\natural}_N)
   = \Det(F^{\frak{o}(S_0)} \,;\, u 1_N + \tilde{\natural}_N).
$$
We can easily check that this $C^{\prime \frak{o}(S_0)}(u)$ is central 
in $U(\frak{o}(S_0))$ for any $u \in \Bbb{C}$.
On the other hand, it is not so easy to calculate its eigenvalue.
However this was given 
through a hard and complicated calculation:

\proclaim {Theorem~3.7 (\cite{I1})}
   \sl We have
   $$
      \pi^{\frak{o}(S_0)}_{\lambda}(C^{\prime \frak{o}(S_0)}(u))
      = \cases
        (u^2 - l_1^2)(u^2 - l_2^2) \cdots (u^2 - l_{N/2}^2), 
        & \qquad \text{$N$: even}, \\
        u (u^2 - l_1^2)(u^2 - l_2^2) \cdots (u^2 - l_{[N/2]}^2), 
        & \qquad \text{$N$: odd}.
        \endcases
   $$
\endproclaim

Comparing this with Theorem~3.6, we have
$C^{\frak{o}(S_0)}(u) = C^{\prime \frak{o}(S_0)}(u)$
(recall that any central element in the universal enveloping algebras
of semisimple Lie algebras is determined by its eigenvalue):

\proclaim {Theorem~3.8 (\cite{W})}
   \sl We have
   $$
      \det(F^{\frak{o}(S_0)} + u \boldkey{1} 
             + \diag\tilde{\natural}_N)
      = \Det(F^{\frak{o}(S_0)} + u \boldkey{1} \,;\, \tilde{\natural}_N).
   $$
\endproclaim

This equality was first shown by A. Wachi in this way.
Namely this proof depends on the two non-trivial results
Theorems~3.5 and~3.7.

However, we can also prove Theorem~3.8 directly 
not using Thereoms~3.5 and~3.7 
(see \cite{I4};
this is similar to the proof of the main theorem in this paper,
but easier).
Conversely, Theorems~3.5 and~3.7 follow
from this Theorem~3.8 immediately.

\vskip .10in

\noindent
{\bf 3.4.}
These results can be generalized in terms of minors:

\proclaim {Theorem~3.9 (\cite{W})} {\bf }
   \sl The following element is central in $U(\frak{o}(S_0))$ 
   for any $u \in \Bbb{C}$:
   $$
      C^{\frak{o}(S_0)}_k(u)
      = \sum_{1 \leq \alpha_1 < \cdots < \alpha_k \leq N}
        \det(\widetilde{F}^{\frak{o}(S_0)}_{\alpha} + u \boldkey{1} 
        + \diag(\tfrac{k}{2}-1,\tfrac{k}{2}-2,\ldots,-\tfrac{k}{2})).
   $$
   Here we put
   $$
      \widetilde{F}^{\frak{o}(S_0)} = 
      \cases
         F^{\frak{o}(S_0)} + \diag(0,\ldots,0,1,\ldots,1), 
         & \qquad \text{$N$: even}, \\
         F^{\frak{o}(S_0)} + \diag(0,\ldots,0,\tfrac{1}{2},1,\ldots,1), 
         & \qquad \text{$N$: odd},
      \endcases
   $$
   where the numbers of $0$'s and $1$'s are equal to $[N/2]$.
\endproclaim

This central element can be rewritten in terms of the symmetrized determinant:

\proclaim {Theorem~3.10 (\cite{W})}
   \sl We have
   $$
      \sum_{1 \leq \alpha_1 < \cdots < \alpha_k \leq N}
      \det(\widetilde{F}^{\frak{o}(S_0)}_{\alpha} + u \boldkey{1} 
             + \diag(\tfrac{k}{2}-1,\tfrac{k}{2}-2,\ldots,-\tfrac{k}{2}))
      = \Det_k(F^{\frak{o}(S_0)} + u \boldkey{1} \,;\, \tilde{\natural}_k).
   $$
\endproclaim

These can be deduced from Theorem~3.8.
See \cite{W} for the details.

\vskip .10in

\noindent
{\it Remarks. }
(1) We can also express $C^{\frak{o}(S_0)}_k(u)$ as
$$
   C^{\frak{o}(S_0)}_k(u)
   = \sum_{1 \leq \alpha_1 < \cdots < \alpha_k \leq N}
     \det(\widehat{F}^{\frak{o}(S_0)}_{\alpha} + u \boldkey{1} 
     + \diag(\tfrac{k}{2},\tfrac{k}{2}-1,\ldots,-\tfrac{k}{2}+1)).
$$
Here $\widehat{F}^{\frak{o}(S_0)}$ is defined by
$$
   \widehat{F}^{\frak{o}(S_0)} = 
   \cases
      F^{\frak{o}(S_0)} - \diag(1,\ldots,1,0,\ldots,0), 
      & \qquad \text{$N$: even}, \\
      F^{\frak{o}(S_0)} - \diag(1,\ldots,1,\tfrac{1}{2},0,\ldots,0), 
      & \qquad \text{$N$: odd},
   \endcases
$$
where the numbers of $1$'s and $0$'s are equal to $[N/2]$.

\vskip .10in

\noindent
(2)
The element $C^{\frak{o}(S_0)}(u)$ is also equal to 
the central element given in \cite{M} in terms of the Sklyanin determinant.
This is seen by comparing their eigenvalues.
See \cite{M}, \cite{MN}, \cite{MNO},
\cite{IU}, \cite{I1}, \cite{W} for the details.

\vskip .10in

\noindent
(3)
The following relation holds for general $S$ \cite{IU}:
$$
   \Det_{2k}(F^{\frak{o}(S)} \,;\, \tilde{\natural}_{2k})
   = \sum_{1 \leq \alpha_1 < \cdots < \alpha_{2k} \leq N}
      \Pf(F^{\frak{o}(S)}S)_{\alpha} \Pf(S^{-1}F^{\frak{o}(S)})_{\alpha}.
\tag 3.1
$$
Here we define the Pfaffian $\Pf Z$ for an alternating matrix $Z = (Z_{ij})$
of size $2k$ by
$$
   \Pf Z = \frac{1}{2^k k!} \sum_{\sigma\in\frak{S}_{2k}} 
   \sgn(\sigma)
   Z_{\sigma(1)\sigma(2)} Z_{\sigma(3)\sigma(4)} 
   \cdots Z_{\sigma(2k-1)\sigma(2k)}.
$$

\vskip .20in

\subhead 4. The case of the symplectic Lie algebras 
\endsubhead
In this section, we introduce the main object of this paper,
namely an analogue of the Capelli determinant 
for the symplectic Lie algebra $\frak{sp}_N$.
We can regard this as the direct counterpart of the element 
$C^{\frak{o}(S_0)}_k(u)$ due to A. Wachi,
but this element is given in terms of the column-permanent
not in terms of the column-determinant.

\vskip .10in 

\noindent
{\bf 4.1.}
First we see the general realization of $\frak{sp}_N$.
Let $J \in \Mat_N(\Bbb{C})$ be a non-degenerate alternating matrix 
of size $N$ (hence $N$ must be even; let us put $n = N/2$).
We can realize the symplectic Lie group
as the isometry group with respect to the bilinear form
determined by $J$:
$$
   Sp(J) = \{ g \in GL_N \,|\, {}^t\!g J g = J \}.
$$
The corresponding Lie algebra is expressed as
$$
   \frak{sp}(J) = \{ Z \in \frak{gl}_N \,|\, {}^t\!Z J + J Z = 0 \}.
$$
As generators of this $\frak{sp}(J)$,
we can take $F^{\frak{sp}(J)}_{ij} = E_{ij} - J^{-1} E_{ji} J$.
We consider the $N \times N$ matrix 
$F^{\frak{sp}(S)} = (F^{\frak{sp}(J)}_{ij})_{1 \leq i,j \leq N}$ 
whose $(i,j)$th entry is this generator $F^{\frak{sp}(J)}_{ij}$.
By a direct calculation, 
this $F^{\frak{sp}(J)}$ satisfies the following relation:

\proclaim {Lemma~4.1} {\bf }
   \sl For any $g \in Sp(J)$, we have
   $$
      \Ad(g) F^{\frak{o}(J)} 
      = {}^t\!g \cdot F^{\frak{sp}(J)} \cdot {}^t\!g^{-1}.
   $$
   Here $\Ad(g) F^{\frak{sp}(J)}$ means the matrix 
   $(\Ad(g) F^{\frak{sp}(J)}_{ij})_{1 \leq i,j \leq N}$.
\endproclaim

Combining this with Proposition~1.9, we have the following proposition:

\proclaim {Proposition~4.2} {\bf }
   \sl The two elements 
   $$
      \Det_k(F^{\frak{sp}(J)} \,;\, a_1,\ldots,a_k), \qquad
      \Per_k(F^{\frak{sp}(J)} \,;\, a_1,\ldots,a_k)
   $$
   are invariant 
   under the adjoint action of $Sp(J)$,
   and in particular this is central in $U(\frak{sp}(J))$.
\endproclaim

Thus, the symmetrized determinant and the symmetrized permanent 
are useful to construct central elements of $U(\frak{sp}(J))$
as in the case of $\frak{gl}_N$.
However, unfortunately,
the column-determinant and the column-permanent do not seem 
so useful for this purpose at least for general $J$.

\vskip .10in


\noindent
{\bf 4.2.}
Let us consider the split realization of the symplectic Lie algebra.
Namely we consider the case
$$
   J = J_0 =
   \left(
   \smallmatrix
      & & & & & 1 \\
      & & & & \lower.6ex\hbox{\iddots} & \\
      & & & 1 & & \\
      & & -1 & & & \\
      & \lower.6ex\hbox{\iddots} & & & & \\
      -1 & & & & &
   \endsmallmatrix
   \right).
$$
It is convenient to introduce the symbols
$$
   i' = N+1-i, \qquad
   \varepsilon(i) = 
   \cases
      -1, & 1 \leq i \leq n, \\
      +1, & n+1 \leq i \leq N,
   \endcases
$$
so that $J_0 = (\varepsilon(j)\delta_{ij'})_{1 \leq i,j \leq N}$ and
$F^{\frak{sp}(J_0)}_{ij} = E_{ij} - \varepsilon(i)\varepsilon(j)E_{j'i'}$.
Moreover the commutation relation of $F^{\frak{sp}(J_0)}_{ij}$ 
is given by
$$\split
   [F^{\frak{sp}(J_0)}_{ij}, F^{\frak{sp}(J_0)}_{kl}] 
   & = F^{\frak{sp}(J_0)}_{il} \delta_{kj} 
     - F^{\frak{sp}(J_0)}_{kj} \delta_{il} \\
   & \qquad
     + \varepsilon(k) \varepsilon(l) F^{\frak{sp}(J_0)}_{l'j} \delta_{ik'}
     + \varepsilon(i) \varepsilon(j) F^{\frak{sp}(J_0)}_{ki'} \delta_{j'l}.
\endsplit
\tag 4.1
$$

This realization $\frak{sp}(J_0)$ is important in the representation theory.
Indeed, we can take a triangular decomposition of $\frak{sp}(J_0)$
simply as follows:
$$
   \frak{sp}(J_0) = \frak{n}^- \oplus \frak{h} \oplus \frak{n}^+.
\tag 4.2
$$
Here $\frak{n}^-$, $\frak{h}$, and $\frak{n}^+$
are the subalgebra of $\frak{sp}(J_0)$ spanned by the elements 
$F^{\frak{sp}(J_0)}_{ij}$ such that
$i>j$, $i=j$, and $i<j$ respectively.
We call this $\frak{sp}(J_0)$ be the ``split realization''
of the symplectic Lie algebra.

The main object of this paper is the following element
of $U(\frak{sp}(J_0))$:
$$
   D^{\frak{sp}(J_0)}_k(u) 
   = \sum_{1 \leq \alpha_1 \leq \cdots \leq \alpha_k \leq N} 
     \frac{1}{\alpha!} \per(\widetilde{F}^{\frak{sp}(J_0)}_{\alpha} 
   + u \boldkey{1}_{\alpha} 
   - \boldkey{1}_{\alpha} 
     \diag(\tfrac{k}{2}-1,\tfrac{k}{2}-2,\ldots,-\tfrac{k}{2})).
$$
Here $\widetilde{F}^{\frak{sp}(J_0)}$ means the matrix 
$$
   \widetilde{F}^{\frak{sp}(J_0)} 
   = F^{\frak{sp}(J_0)} - \diag(0,\ldots,0,1,\ldots,1),
$$
where the numbers of $0$'s and $1$'s are equal to $n$.

\proclaim {Theorem~4.3} {\bf }
   \sl The element $D^{\frak{sp}(J_0)}_k(u)$ is central in $U(\frak{sp}(J_0))$
   for any $u \in \Bbb{C}$.
\endproclaim

The eigenvalue of $D^{\frak{sp}(J_0)}_k(u)$ 
on the irreducible representations of $\frak{sp}(J_0)$ 
can be calculated easily 
by noting the triangular decomposition (4.2):

\proclaim {Theorem~4.4} {\bf }
   \sl For the representation $\pi^{\frak{sp}(J_0)}_{\lambda}$ 
   of $\frak{sp}(J_0)$ determined by the partition
   $\lambda = (\lambda_1,\ldots,\lambda_n)$,
   we have
   $$\align
      & \pi^{\frak{sp}(J_0)}_{\lambda}(D^{\frak{sp}(J_0)}_k(u)) 
      = \sum_{l=0}^k 
        \sum \Sb 
                 1 \leq \alpha_1 \leq \cdots \leq \alpha_l \leq n \\
                 n+1 \leq \alpha_{l+1} \leq \cdots \leq \alpha_k \leq N
             \endSb \\
      & \qquad 
     (u + \lambda_{\alpha_1} - \tfrac{k}{2} +1)
     (u + \lambda_{\alpha_2} - \tfrac{k}{2} +2)
     \cdots
     (u + \lambda_{\alpha_l} - \tfrac{k}{2} +l) \\
     & \qquad
     \cdot
     (u - \lambda_{\alpha'_{l+1}} - \tfrac{k}{2} +l)
     (u - \lambda_{\alpha'_{l+2}} - \tfrac{k}{2} +l+1)
     \cdots
     (u - \lambda_{\alpha'_k} + \tfrac{k}{2} -1).
   \endalign
   $$
\endproclaim

To prove Theorem~4.3, we additionally consider the following element:
$$
   D^{\prime\frak{sp}(J_0)}_k(u) 
   = \Per_k (F^{\frak{sp}(J_0)} + u \boldkey{1} \,;\, \tilde{\natural}_k)
   = \Per_k (F^{\frak{sp}(J_0)} \,;\, u 1_N + \tilde{\natural}_k).
$$
It is obvious from Proposition~4.2 that 
this $D^{\prime\frak{sp}(J_0)}_k(u)$ 
is central in $U(\frak{sp}(J_0))$ for any $u \in \Bbb{C}$.
However it is not so easy to calculate the eigenvalue 
of $D^{\prime\frak{sp}(J_0)}_k(u)$ directly.

Actually these two elements  $D^{\frak{sp}(J_0)}_k(u)$ and 
$D^{\prime\frak{sp}(J_0)}_k(u)$ are equal:

\proclaim {Theorem~4.5} {\bf }
   \sl We have $D^{\frak{sp}(J_0)}_k(u) = D^{\prime\frak{sp}(J_0)}_k(u)$,
   namely
   $$\align
     \sum_{1 \leq \alpha_1 \leq \cdots \leq \alpha_k \leq N}
     \frac{1}{\alpha!} 
     \per(\widetilde{F}^{\frak{sp}(J_0)}_{\alpha} + u \boldkey{1}_{\alpha} 
        - \boldkey{1}_{\alpha} 
          & \diag(\tfrac{k}{2}-1, \tfrac{k}{2}-2, \ldots, -\tfrac{k}{2})) \\
     & \qquad 
    = \Per_k (F^{\frak{sp}(J_0)} + u \boldkey{1} \,;\, \tilde{\natural}_k).
   \endalign
   $$
\endproclaim

Theorem~4.3 is immediate from this. 
Moreover, using Theorems~4.4 and 4.5, 
we can easily see the eigenvalue of $D^{\prime\frak{sp}(J_0)}_k(u)$.
Thus, as in the case of $\frak{o}(S_0)$,
this Theorem~4.5 settles two problems at the same time.

The remainder of this paper is devoted to the proof of Theorem~4.5.

\vskip .10in

\noindent
{\it Remarks.}
(1) We can also express $D^{\frak{sp}(J_0)}_k(u)$ as 
$$
   D^{\frak{sp}(J_0)}_k(u) 
   = \sum_{1 \leq \alpha_1 \leq \cdots \leq \alpha_k \leq N} 
     \frac{1}{\alpha!} \per(\widehat{F}^{\frak{sp}(J_0)}_{\alpha} 
   + u \boldkey{1}_{\alpha} 
   - \boldkey{1}_{\alpha} 
     \diag(\tfrac{k}{2},\tfrac{k}{2}-1,\ldots,-\tfrac{k}{2}+1)).
$$
Here we put
$\widehat{F}^{\frak{sp}(J_0)} 
= F^{\frak{sp}(J_0)} + \diag(1,\ldots,1,0,\ldots,0)$.

\vskip .10in

\noindent
(2) 
Considering the generating function of $D^{\frak{sp}(J_0)}_k(u)$,
we can rewrite Theorem~4.4 more simply (see \cite{I3}).
Moreover the right hand side of Theorem~4.4 can be regarded as
the complete symmetric polynomials associated to a kind of factorial power.
See \cite{I5} for the details.

\vskip .10in

\noindent
(3)
From Theorem~4.4, we see that $D^{\frak{sp}(J_0)}_k(0)$ is 
equal to $D_k$ defined in \cite{MN}.
Moreover, as the counterpart of (3.1),
we have the relation 
$$
   \Per_{2k}(F^{\frak{sp}(J)} \,;\, \tilde{\natural}_{2k})
   = \sum_{1 \leq \alpha_1 \leq \cdots \leq \alpha_{2k} \leq N}
     \frac{1}{\alpha!}
      \Hf(F^{\frak{sp}(J)}J)_{\alpha} \Hf(J^{-1}F^{\frak{sp}(J)})_{\alpha}
$$
for general $J$.
Here we define the Hafnian $\Hf Z$ for a symmetric matrix $Z = (Z_{ij})$
of size $2k$ by
$$
   \Hf Z = \frac{1}{2^k k!} \sum_{\sigma\in\frak{S}_{2k}} 
   Z_{\sigma(1)\sigma(2)} Z_{\sigma(3)\sigma(4)} 
   \cdots Z_{\sigma(2k-1)\sigma(2k)}.
$$
This is deduced from Theorem~5.1 in \cite{MN}.

\vskip .20in

\subhead 5. Proof of the main theorem
\endsubhead
Let us show Theorem~4.5
using the symmetric tensor algebra.
This proof is similar to that of Theorem~1.8,
but more complicated.
Namely we need the variable transformation method
developed in \cite{IU}, \cite{I2}, \cite{I3}.

Hereafter, we omit the superscript $\frak{sp}(J_0)$.
Namely we denote $F^{\frak{sp}(J_0)}$,
$F^{\frak{sp}(J_0)}_{ij}$,
$D^{\frak{sp}(J_0)}_k(u)$,
and $D^{\prime\frak{sp}(J_0)}_k(u)$ simply by
$F$, $F_{ij}$, $D_k(u)$, and $D'_k(u)$ respectively.

\vskip .10in


\noindent
{\bf 5.1.}
First, let us express both sides of Theorem~4.5 
using the symmetric tensor algebra
$S_{2N} = S(\Bbb{C}^N \oplus \Bbb{C}^N)$.
Let $e_1,\ldots,e_N, e^*_1,\ldots,e^*_N$
be the standard generators of $S_{2N}$.

In the extended algebra $S_{2N} \otimes U(\frak{sp}(J_0))$, we put
$$
   \varXi = \sum_{i,j=1}^N e_i e^*_j F_{ij}, \qquad
   \varXi(u) = \sum_{i,j=1}^N e_i e^*_j F_{ij}(u), \qquad
   \tau = \sum_{i=1}^N e_i e^*_i,
$$
so that $\varXi(u) = \varXi + u \tau$.
Then, by (2.7), we can express $D'_k(u)$ as 
$$\split
   D'_k(u) 
   & = \Per_k(F + u \boldkey{1} \,;\, 
     -\tfrac{k}{2}+1, -\tfrac{k}{2}+2, \ldots, \tfrac{k}{2}-1, 0) \\
   & = \frac{1}{k!} 
   \big< 
      \varXi(u-\tfrac{k}{2}+1)
      \varXi(u-\tfrac{k}{2}+2) 
      \cdots
      \varXi(u+\tfrac{k}{2}-1) \cdot \varXi(u) 
   \big> \\
   & = \frac{1}{k!} 
   \big< 
      \varXi^{\overline{k-1}}(u-\tfrac{k}{2}+1) \cdot \varXi(u) 
   \big>.
\endsplit
\tag 5.1
$$
Here $\varXi^{\overline{k}}(u)$ means 
the ``rising factorial power''
$$
   \varXi^{\overline{k}}(u) = \varXi(u) \varXi(u+1) \cdots \varXi(u+k-1).
$$

Let us express $D_k(u)$ similarly.
We put $\eta_j(u) = \sum_{i=1}^N e_i F_{ij}(u)$
and  $\eta^{\dagger}_j(u) = \eta_j(u) e^*_j$,
so that $\sum_{j=1}^N \eta^{\dagger}_j(u) = \varXi(u)$.
Moreover, 
we put $\tilde{\eta}_j(u) = \sum_{i=1}^N e_i \widetilde{F}_{ij}(u)$
and  $\tilde{\eta}^{\dagger}_j(u) = \tilde{\eta}_j(u) e^*_j$.
Then $\eta^{\dagger}_j(u)$ and $\tilde{\eta}^{\dagger}_j(u)$
are related by
$$
   \tilde{\eta}^{\dagger}_i (u) = 
   \cases
      \eta^{\dagger}_i (u), \quad &  \phantom{n+{}} 1 \leq i \leq n, \\
      \eta^{\dagger}_i (u-1), \quad & n+1 \leq  i \leq N.
   \endcases
\tag 5.2
$$
In this notation, we have
$$\align
   & \per(\widetilde{F}_{\alpha} + u \boldkey{1}_{\alpha}
   - \boldkey{1}_{\alpha} 
   \diag(\tfrac{k}{2}-1,\tfrac{k}{2}-2,\ldots,-\tfrac{k}{2})) \\
   & \qquad = 
   \big<
         \tilde{\eta}^{\dagger}_{\alpha_1}(u-\tfrac{k}{2}+1)
         \tilde{\eta}^{\dagger}_{\alpha_2}(u-\tfrac{k}{2}+2) 
         \cdots
         \tilde{\eta}^{\dagger}_{\alpha_k}(u+\tfrac{k}{2})
   \big>
\endalign
$$
by (2.8).
Thus we can express $D_k(u)$ as
$$
   D_k(u) 
   = \sum_{1 \leq \alpha_1 \leq \cdots \leq \alpha_k \leq N}
   \frac{1}{\alpha!}
   \big< 
         \tilde{\eta}^{\dagger}_{\alpha_1}(u-\tfrac{k}{2}+1)
         \tilde{\eta}^{\dagger}_{\alpha_2}(u-\tfrac{k}{2}+2) 
         \cdots
         \tilde{\eta}^{\dagger}_{\alpha_k}(u+\tfrac{k}{2}) 
   \big>.
\tag 5.3
$$

\vskip .10in

\noindent
{\it Remark.}
Recall that
$\varXi(u)$ and $\varXi(w)$ are commutative for any $u$, $w \in \Bbb{C}$,
because $\tau$ is central.

\vskip .10in


\noindent
{\bf 5.2.}
By (5.1) and (5.3), 
our goal $D_k(u) = D'_k(u)$ can be expressed as
$$
   \sum_{1 \leq \alpha_1 \leq \cdots \leq \alpha_k \leq N}
   \frac{k!}{\alpha!}
   \big< 
         \tilde{\eta}^{\dagger}_{\alpha_1}(u-\tfrac{k}{2}+1)
         \tilde{\eta}^{\dagger}_{\alpha_2}(u-\tfrac{k}{2}+2) 
         \cdots
         \tilde{\eta}^{\dagger}_{\alpha_k}(u+\tfrac{k}{2}) 
   \big>
   =
   \big< 
      \varXi^{\overline{k-1}}(u-\tfrac{k}{2}+1) \cdot \varXi(u) 
   \big>.
$$
Replacing $u$ by $u+\frac{k}{2}-1$, 
we can rewrite this simply as
$$
   \sum_{1 \leq \alpha_1 \leq \cdots \leq \alpha_k \leq N}
   \frac{k!}{\alpha!}
   \big< 
         \tilde{\eta}^{\dagger}_{\alpha_1}(u)
         \tilde{\eta}^{\dagger}_{\alpha_2}(u+1) 
         \cdots
         \tilde{\eta}^{\dagger}_{\alpha_k}(u+k-1) 
   \big>
   =
   \big< 
      \varXi^{\overline{k-1}}(u) \cdot \varXi(u+\tfrac{k}{2}-1) 
   \big>.
$$
Let us prove Theorem~4.5 in this form.
Namely, we hereafter aim the following relation:

\proclaim{Lemma~5.1} \sl
   We have
   $$
      \big< W_k(u) \big> = \big< W'_k(u) \big>.
   $$
   Here we put
   $$\align
      W_k(u) & = 
      \sum_{1 \leq \alpha_1 \leq \cdots \leq \alpha_k \leq N}  
      \frac{k!}{\alpha!}
            \tilde{\eta}^{\dagger}_{\alpha_1}(u)
            \tilde{\eta}^{\dagger}_{\alpha_2}(u+1) 
            \cdots
            \tilde{\eta}^{\dagger}_{\alpha_k}(u+k-1), \\
      W'_k(u) 
      & = \varXi^{\overline{k-1}}(u) \cdot \varXi(u+\tfrac{k}{2}-1)
      = \varXi^{\overline{k}}(u) 
        - \frac{k}{2} \varXi^{\overline{k-1}}(u) \tau.
   \endalign
   $$
\endproclaim

\vskip .10in


\noindent
{\bf 5.3.}
In Sections~5.3--5.5, we will study the relation
between $W_k(u)$ and the factorial powers of $\varXi(u)$.
First, 
$W_k(u)$ is expressed in terms of $\eta^{\dagger}_{i}(u)$ as
$$\split
  W_k(u)
   & = \sum_{l=0}^k 
   \sum \Sb 1 \leq \alpha_1 \leq \cdots \leq \alpha_l \leq n \\
            n+1 \leq \alpha_{l+1} \leq \cdots \leq \alpha_k \leq N \endSb
   \frac{k!}{\alpha!} 
         \eta^{\dagger}_{\alpha_1}(u)
         \eta^{\dagger}_{\alpha_2}(u+1) 
         \cdots
         \eta^{\dagger}_{\alpha_l}(u+l-1) \\
         & \qquad\qquad 
         \cdot
         \eta^{\dagger}_{\alpha_{l+1}}(u+l-1)
         \eta^{\dagger}_{\alpha_{l+2}}(u+l) 
         \cdots
         \eta^{\dagger}_{\alpha_k}(u+k-2).
\endsplit
\tag 5.4
$$
Note that $\eta_{i}(u)$ and $\eta^{\dagger}_{i}(u)$ 
satisfy the following commutation relation.
This is deduced from~(4.1) by a direct calculation.

\proclaim{Lemma~5.2}
   \sl We have
   $$\align
      \eta_j(u) \eta_l(u+1) - \eta_l(u) \eta_j(u+1)
      & = \varTheta  J_0^{jl}
        = \varepsilon(j) \varTheta \delta_{j'l}, \\
      \eta^{\dagger}_j(u) \eta^{\dagger}_l(u+1) 
      - \eta^{\dagger}_l(u) \eta^{\dagger}_j(u+1)
      & = \varTheta  J_0^{jl} e^*_j e^*_l
        = \varepsilon(j) \varTheta \delta_{j'l} e^*_j e^*_l. 
   \endalign
   $$
   Here $J_0^{ij}$ means the $(i,j)$th entry of the matrix $J_0^{-1}$,
   and $\varTheta$ is defined by
   $$
      \varTheta 
      = \sum_{a,b=1}^N \varepsilon(b) e_a e_b F_{ab'}.
   $$
\endproclaim

\noindent
{\it Remark.}
For any $u \in \Bbb{}C$, we have 
$$
   \varTheta = \sum_{b=1}^N \varepsilon(b) \eta_b(u) e_b.
$$

\proclaim{Corollary~5.3}
   \sl When $1 \leq i,j \leq n$ or $n+1 \leq i,j \leq N$, 
   we have
   $$
  \eta^{\dagger}_i(u) \eta^{\dagger}_j(u+1) 
   = \eta^{\dagger}_j(u) \eta^{\dagger}_i(u+1).
   $$
\endproclaim

Noting this relation, 
we consider the two elements
$$\align
   \varXi_{-}(u) & = \sum_{j=1}^n \eta^{\dagger}_j(u)
   = \sum_{i=1}^N \sum_{j=1}^n e_i e^*_j F_{ij}(u), \\
   \varXi_{+}(u) & = \sum_{j=n+1}^N \eta^{\dagger}_j(u)
   = \sum_{i=1}^N \sum_{j=n+1}^N e_i e^*_j F_{ij}(u).
\endalign
$$
Then we have $\varXi_{-}(u) + \varXi_{+}(u) = \varXi(u)$.
Moreover we put
$$\align
   \varXi^{\overline{k}}_{-}(u)
   & = \varXi_{-}(u) \varXi_{-}(u+1) \cdots \varXi_{-}(u+k-1), \\
   \varXi^{\overline{k}}_{+}(u)
   & = \varXi_{+}(u) \varXi_{+}(u+1) \cdots \varXi_{+}(u+k-1).
\endalign
$$
By Corollary~5.3, these factorial powers can be expanded as 
$$\align
   \varXi^{\overline{k}}_{-}(u)
   & = \sum_{1 \leq \alpha_1 \leq \cdots \leq \alpha_k \leq n}
     \frac{k!}{\alpha!}
     \eta^{\dagger}_{\alpha_1}(u) 
     \eta^{\dagger}_{\alpha_2}(u+1)
     \cdots 
     \eta^{\dagger}_{\alpha_k}(u+k-1), \\
   \varXi^{\overline{k}}_{+}(u)
   & = \sum_{n+1 \leq \beta_1 \leq \cdots \leq \beta_k \leq N}
     \frac{k!}{\beta!}
     \eta^{\dagger}_{\beta_1}(u) 
     \eta^{\dagger}_{\beta_2}(u+1)
     \cdots 
     \eta^{\dagger}_{\beta_k}(u+k-1).
\endalign
$$
Thus we can rewrite (5.4) simply as 
$$
   W_k (u) = \sum_{l \geq 0} {k \choose l} 
   \varXi_{-}^{\overline{l}}(u)
   \varXi_{+}^{\overline{k-l}}(u+l-1).
$$

\vskip .10in


\noindent
{\bf 5.4.}
Let us consider an analogue of $W_k (u)$:
$$
   V_k (u) = \sum_{l \geq 0} {k \choose l} 
   \varXi_{-}^{\overline{l}}(u)
   \varXi_{+}^{\overline{k-l}}(u+l).
$$
This $V_k (u)$ is related to $W_k (u)$ as follows:

\proclaim{Lemma~5.4}\sl
   We have
   $$
      V_k(u) = W_k(u) + k V_{k-1}(u) \tau_{+}.
   $$
   Here we define $\tau_{-}$ and $\tau_{+}$ by
   $$
      \tau_{-} = \sum_{i=1}^n e_i e^*_i, \qquad
      \tau_{+} = \sum_{i=n+1}^{N} e_i e^*_i,
   $$
   so that $\tau = \tau_{-} + \tau_{+}$,
   $\varXi_{-}(u) = \varXi_{-}(0) + u \tau_{-}$, and 
   $\varXi_{+}(u) = \varXi_{+}(0) + u \tau_{+}$.
\endproclaim

\vskip .10in

\noindent
{\it Remark.}
These $\tau_{-}$ and $\tau_{+}$ are obviously central. 
Hence $\varXi_{-}(u)$ and $\varXi_{-}(w)$ commute with one another
for any $u$, $w \in \Bbb{C}$.
Similarly $\varXi_{+}(u)$ and $\varXi_{+}(w)$ are commutative.

\vskip .10in

\noindent
{\it Proof of Lemma~{\sl 5.4.}}
Since $\varXi_{+}(u) = \varXi_{+}(0) + u \tau_{+}$,
we have
$$
   \varXi_{+}^{\overline{k}}(u) - \varXi_{+}^{\overline{k}}(u-1)
   = k \varXi_{+}^{\overline{k-1}}(u) \tau_{+}.
$$
Hence, we have
\TagsOnRight
$$\align
   V_k(u) - W_k(u)
   & = \sum_{l \geq 0} {k \choose l} \varXi_{-}^{\overline{l}} (u)
   \cdot
   \{ \varXi_{+}^{\overline{k-l}}(u+l) 
   - \varXi_{+}^{\overline{k-l}}(u+l-1) \} \\
   & = \sum_{l \geq 0} {k \choose l} \varXi_{-}^{\overline{l}} (u)
   \cdot (k-l) \varXi_{+}^{\overline{k-l-1}}(u+l) \tau_{+}\\
   & = k \sum_{l \geq 0} {k-1 \choose l} \varXi_{-}^{\overline{l}} (u)
   \varXi_{+}^{\overline{k-l-1}}(u+l) \tau_{+} \\
   & = k V_{k-1}(u) \tau_{+}.
\tag"\qed"
\endalign
$$
\TagsOnLeft

\vskip .10in


\noindent
{\bf 5.5.}
Let us study the relation between $V_k(u)$ 
and the factorial powers of $\varXi(u)$.
First, the following commutation relations
are easy from Lemma~5.2: 

\proclaim{Lemma~5.5}\sl
   We have
   $$
      \varXi_{+}(u-1) \varXi_{-}(u) -  \varXi_{-}(u-1) \varXi_{+}(u)
      = \varTheta \rho^*. 
   $$
   Here we put $\rho^* = \sum_{i=1}^n e^*_i e^*_{i'}$.
\endproclaim

\proclaim{Lemma~5.6}\sl
   We have
   $$
      \eta_j(u) \varTheta = \varTheta \eta_j(u+2),
   $$
   and in particular
   $$
      \varXi(u) \varTheta = \varTheta \varXi(u+2), \qquad
      \varXi_{-}(u) \varTheta = \varTheta \varXi_{-}(u+2), \qquad
      \varXi_{+}(u) \varTheta = \varTheta \varXi_{+}(u+2).
   $$
\endproclaim

Note that $\rho^*$ is central in $S_{2N} \otimes U(\frak{sp}(J_0))$.
Thus, the following relation is obtained from Lemmas~5.5 and 5.6
by a simple calculation.

\proclaim{Lemma~5.7}\sl
   We have
   $$
      \varXi_{+}^{\overline{k}}(u) \varXi_{-}(u+k)
        - \varXi_{-}(u) \varXi_{+}^{\overline{k}}(u+1)
      = k \varXi_{+}^{\overline{k-1}}(u) \varTheta \rho^*.
   $$
\endproclaim

Moreover, we have the following relation:

\proclaim{Lemma~5.8}\sl
   We have
   $$
      V_k(u) \varXi(u+k) - V_{k+1}(u) 
      = k V_{k-1}(u) \varTheta \rho^*.
   $$
\endproclaim

\vskip .10in

\noindent
{\it Proof of Lemma~{\sl 5.8.}}
First, we have
$$\align
   & V_k(u) \varXi_{-}(u+k) - \varXi_{-}(u) V_k(u+1) \\
\allowdisplaybreak
   & \qquad
   = \sum_{l \geq 0} {k \choose l} 
   \varXi_{-}^{\overline{l}}(u) 
   \varXi_{+}^{\overline{k-l}}(u+l) \varXi_{-}(u+k) 
   - \sum_{l \geq 0} {k \choose l} 
   \varXi_{-}^{\overline{l+1}}(u) 
   \varXi_{+}^{\overline{k-l}}(u+l+1) \\
\allowdisplaybreak
   & \qquad
   = \sum_{l \geq 0} {k \choose l} 
   \varXi_{-}^{\overline{l}}(u) \cdot
   \{ \varXi_{+}^{\overline{k-l}}(u+l) \varXi_{-}(u+k) 
   - \varXi_{-}(u+l) \varXi_{+}^{\overline{k-l}}(u+l+1) \} \\
\allowdisplaybreak
   & \qquad
   = \sum_{l \geq 0} {k \choose l} 
   \varXi_{-}^{\overline{l}}(u) \cdot
   (k-l) \varXi_{+}^{\overline{k-l-1}}(u+l) \varTheta \rho^* \\
\allowdisplaybreak
   & \qquad
   = \sum_{l \geq 0} k {k-1 \choose l} 
   \varXi_{-}^{\overline{l}}(u) 
   \varXi_{+}^{\overline{k-l-1}}(u+l) \varTheta \rho^* \\
\allowdisplaybreak
   & \qquad
   = k V_{k-1}(u) \varTheta \rho^*.
\endalign
$$
Here we used Lemma~5.7 for the third equality. 
Moreover, we have
$$\align
   & V_k(u) \varXi_{+}(u+k) + \varXi_{-}(u) V_k(u+1) \\
\allowdisplaybreak
   & \quad= \sum_{l \geq 0} {k \choose l} 
   \varXi_{-}^{\overline{l}}(u) \varXi_{+}^{\overline{k-l+1}}(u+l)
   + \sum_{l \geq 0} {k \choose l} 
   \varXi_{-}^{\overline{l+1}}(u) 
   \varXi_{+}^{\overline{k-l}}(u+l+1) \\
\allowdisplaybreak
   & \quad= \sum_{l \geq 0} {k \choose l} 
   \varXi_{-}^{\overline{l}}(u) \varXi_{+}^{\overline{k+1-l}}(u+l)
   + \sum_{l \geq 0} {k \choose l-1} 
   \varXi_{-}^{\overline{l}}(u) 
   \varXi_{+}^{\overline{k+1-l}}(u+l) \\
\allowdisplaybreak
   & \quad = \sum_{l \geq 0} {k+1 \choose l} 
   \varXi_{-}^{\overline{l}}(u) 
   \varXi_{+}^{\overline{k+1-l}}(u+l) \\
\allowdisplaybreak
   & \quad = V_{k+1}(u),
\endalign
$$
because ${k \choose l} + {k \choose l-1} = {k+1 \choose l}$.
Thus we have
\TagsOnRight
$$\align
   & V_k(u) \varXi(u+k) - V_{k+1}(u) \\
   & \qquad
   = V_k(u) \{ \varXi_{-}(u+k) + \varXi_{+}(u+k) \}
   - \{ V_k(u) \varXi_{+}(u+k) + \varXi_{-}(u) V_k(u+1) \} \\
   & \qquad
   = V_k(u) \varXi_{-}(u+k) - \varXi_{-}(u) V_k(u+1) \\
   & \qquad
   = k V_{k-1}(u) \varTheta \rho^*.
\tag"\qed"
\endalign
$$
\TagsOnLeft

\vskip .10in

Using this, we can show the following expansion:

\proclaim{Lemma~5.9}\sl
   We have
   $$
      \varXi^{\overline{k}}(u)
      = \sum_{l \geq 0} R^k_l V_{k-2l}(u) \varTheta^l \rho^{* l}.
   $$
   Here we put
   $$
      R^k_l 
      = {k \choose 2l} (2l-1)!!
   $$
   with $(2l-1)!! = (2l-1) (2l-3) \cdots 1$
   {\rm (}we put $(-1)!! = 1$ when $l=0${\rm )}. 
\endproclaim

\vskip .10in

\noindent
{\it Proof of Lemma~{\sl 5.9.}}
This can be proved by induction on $k$.
First, the case $k=0$ is easy.
Next, by assuming the case $k=m$,
the case $k=m+1$ is deduced as follows:
$$\align
   \varXi^{\overline{m+1}}(u)
   & = \varXi^{\overline{m}}(u) \varXi(u+m) \\
\allowdisplaybreak
   & = \sum_{l \geq 0} 
   R^m_l V_{m-2l}(u) \varTheta^l \rho^{*l} \varXi(u+m) \\
\allowdisplaybreak
   & = \sum_{l \geq 0} 
   R^m_l V_{m-2l}(u) \varXi(u+m-2l) \varTheta^l  \rho^{*l}  \\
\allowdisplaybreak
   & = \sum_{l \geq 0} R^m_l V_{m-2l+1}(u) \varTheta^l \rho^{*l}
     + \sum_{l \geq 0}
     (m-2l) R^m_l V_{m-2l-1}(u) \varTheta^{l+1} \rho^{*l+1} \\
\allowdisplaybreak
   & = \sum_{l \geq 0} R^m_l V_{m-2l+1}(u) \varTheta^l \rho^{*l}
     + \sum_{l \geq 1} 
     (m-2l+2) R^m_{l-1} V_{m-2l+1}(u) \varTheta^{l} \rho^{*l} \\
\allowdisplaybreak
   & = \sum_{l \geq 0} R^{m+1}_l V_{m-2l+1}(u) \varTheta^l \rho^{*l}.
\endalign
$$
Here we used Lemma~5.8 for the fourth equality.
To show the last equality,
we also used the relations
$$
   R^{m+1}_0 = R^m_0 = 1, \qquad 
   R^{m+1}_l = R^m_l + (m-2l+2) R^m_{l-1}.
$$
These relations themselves are immediate from the definition of $R^k_l$.
\hfill\qed

\vskip .10in

The coefficient $R^k_l$ also appears in the expansions
$$
   u^{\overline{k}} 
   = \sum_{l \geq 0} (-)^l R^k_l u^{\overline{\overline{k-l}}}, \qquad
   u^{\overline{\overline{k}}} 
   = \sum_{l \geq 0} R^{k+l-1}_l u^{\overline{k-l}}.
\tag 5.5
$$
Here $u^{\overline{k}}$ and $u^{\overline{\overline{k}}}$ 
mean the two factorial powers
$$
   u^{\overline{k}} = u (u+1) (u+2) \cdots (u+k-1), \qquad 
   u^{\overline{\overline{k}}} = u (u+2) (u+4) \cdots (u+2k-2).
$$
In particular, we have
$$
   \sum_{l \geq 0} (-)^l R^k_l R^{k-2l}_{m-l} = \delta_{m,0}.
$$
By noting this, the following is immediate from Lemma~5.9:

\proclaim{Lemma~5.10}\sl
   We have
   $$
      V_k(u) 
      = \sum_{l \geq 0} 
        (-)^l R^k_l \varXi^{\overline{k-2l}}(u) \varTheta^l \rho^{* l}.
   $$
\endproclaim

\vskip .10in

\noindent
{\bf 5.6.}
Next we consider the following relations:

\proclaim {Lemma~5.11} {\bf }
   \sl We have
   $$
      k
      \big< 
         \varXi^{\overline{\overline{k-1}}}(u) 
         \varTheta^{l} \rho^{* l} \tau^m
      \big>
      +l
      \big< 
         \varXi^{\overline{\overline{k}}}(u) 
         \varTheta^{l-1} \rho^{* l-1} \tau^m \omega 
      \big>
      = 0.
   $$
\endproclaim

\proclaim {Lemma~5.12} {\bf }
   \sl We have
   $$
      k \big< W'_{k-1}(u) \varTheta^{l} \rho^{* l} \big>
      = l \big< \varXi^{\overline{k}}(u) 
      \varTheta^{l-1} \rho^{* l-1} \omega \big>.
   $$
\endproclaim

Here $\omega$ is the central element defined by
$$
   \omega = \sum_{i=1}^N \varepsilon(i) e_i e^*_i
   = -\tau_{-} + \tau_{+}.
$$
Moreover $\varXi^{\overline{\overline{k}}}(u)$ means 
the factorial power
$$
   \varXi^{\overline{\overline{k}}}(u)
   = \varXi(u) \varXi(u+2) \varXi(u+4) \cdots \varXi(u+2k-2). 
$$

To prove Lemma~5.11, we use a variable transformation and Lemma~2.1.
We put
$$
   g = 
   \pmatrix
      a \boldkey{1} & b {}^t\! J_0 \\
      c {}^t\! J_0^{-1} & d \boldkey{1}
   \endpmatrix \in GL_{2N},
$$
so that
$$
   {}^t\!g^{-1} = \frac{1}{ad-bc}
   \pmatrix
      d \boldkey{1} & -c J_0^{-1} \\
      -b J_0 & a \boldkey{1}
   \endpmatrix.
$$
These $g$ and ${}^t\!g^{-1}$ naturally act on $\Bbb{C}^{2N}$, 
the space spanned by the formal variables 
$e_1,\ldots, e_N, e^*_1,\ldots, e^*_N$.
Let us consider their extended actions on $S_{2N} = S(\Bbb{C}^{2N})$
and moreover on $S_{2N} \otimes U(\frak{sp}(J_0))$
as automorphisms.
Then, we have the relations
$$
\aligned
   g(\tau) & = (ad-bc) \tau, \\
   g(\rho) & = a^2 \rho + c^2 \rho^* + ac \omega, \\
   g(\varXi) & = (ad+bc) \varXi + ab \varTheta + cd \varTheta^*, \\
   g(\varTheta) & = a^2 \varTheta + c^2 \varTheta^* + 2ac \varXi, 
\endaligned
\quad
\aligned
   g(\omega) & = (ad+bc) \omega + 2ab \rho + 2cd \rho^*, \\
   g(\rho^*) & = b^2 \rho + d^2 \rho^* + bd \omega, \\
   & \\
   g(\varTheta^*) &= b^2 \varTheta + d^2 \varTheta^* + 2bd \varXi.
\endaligned
\tag 5.6
$$
Here we define $\varTheta^*$ and $\rho$ by
$$
   \varTheta^*
   = - \sum_{i,j=1}^N \varepsilon(i) e^*_i e^*_j F_{i'j}, \qquad
   \rho = - \sum_{i=1}^n e_i e_{i'}.
$$
Moreover we have
$$
   {}^t\!g^{-1}(\tau) = (ad-bc)^{-1} \tau.
$$

\vskip .10in

\noindent
{\it Remark.}
To show (5.6),
it is convenient to consider 
the ``row vectors'' $e = (e_1,\ldots,e_N)$ and 
$e^* = (e^*_1,\ldots,e^*_N)$,
so that 
$$\gathered
   \tau = e {}^t\!e^*, \qquad
   \varXi = e F {}^t\!e^*, \qquad
   \varTheta = e F J_0 {}^t\!e, \qquad
   \varTheta^* = e^* {}^t\!J_0^{-1} F {}^t\!e^*, \\
   \omega = e K_0 {}^t\!e^*, \qquad
   \rho = \frac{1}{2} e K_0 J_0 {}^t\!e, \qquad
   \rho^* = \frac{1}{2} e^* {}^t\!J_0^{-1} K_0 {}^t\!e^*.
\endgathered
$$
Here $K_0$ means the matrix $K_0 = \diag(-1,\ldots,-1,1,\ldots,1)$.

\vskip .10in

\noindent
{\it Proof of Lemma~{\sl 5.11.}}
Let us suppose that $a=b=d=1$ and $c=0$. 
Then (5.6) is rewritten as
$$\gathered
   g(\tau) = \tau, \qquad
   g(\omega) = \omega + 2 \rho, \qquad
   g(\rho) = \rho, \qquad
   g(\rho^*) = \rho + \rho^* + \omega, \\
   g(\varXi) = \varXi + \varTheta, \qquad
   g(\varTheta) = \varTheta, \qquad
   g(\varTheta^*) = \varTheta + \varTheta^* + 2 \varXi.
\endgathered
$$
In particular we have
$$
   g(\varXi(u)) = \varXi(u) + \varTheta, \qquad
   g(\varTheta) = \varTheta, \qquad
   g(\omega - 2\rho^*) = - (\omega + 2\rho^*), \qquad
   g(\tau) = \tau.
$$
Using Lemma~5.6, we can show a kind of binomial expansion:
$$\align
   g(\varXi^{\overline{\overline{k}}}(u))
   & = (\varXi(u) + \varTheta) (\varXi(u+2) + \varTheta)
   \cdots
   (\varXi(u + 2k-2) + \varTheta) \\
   & = \sum_{s \geq 0} {k \choose s}
   \varXi(u) \varXi(u+2) \cdots \varXi(u+2s-2)  
   \cdot
   \varTheta^{k-s} \\
   & = \sum_{s \geq 0} {k \choose s}
   \varXi^{\overline{\overline{s}}}(u)
   \cdot
   \varTheta^{k-s}.
\endalign
$$
Moreover, expanding the equality 
$g((\omega - 2\rho^*)^l) = (-)^l (\omega + 2\rho^*)^l$, 
we have
$$
   g \big( \sum_{r \geq 0} {l \choose r} 
   \omega^r (- 2\rho^*)^{l-r} \big)
   = (-)^l \sum_{r \geq 0} {l \choose r} 
   \omega^r (2\rho^*)^{l-r}.
$$
We also have $g(\varTheta^{l-1}) = \varTheta^{l-1}$ 
and $g(\tau^m) = \tau^m$.
Multiplying these equalities, we have
$$\align
   & g \big( \varXi^{\overline{\overline{k}}}(u) \varTheta^{l-1}
   \sum_{r \geq 0} {l \choose r} 
   \omega^r (- 2\rho^*)^{l-r} \tau^m \big) \\
   & \qquad =
   (-)^{l} 
   \sum_{s \geq 0} {k \choose s} \varXi^{\overline{\overline{s}}}(u) 
   \varTheta^{k-s+l-1} 
   \sum_{r \geq 0} {l \choose r} \omega^r (2\rho^*)^{l-r} \tau^m.
\endalign
$$
Take the inner product with $\tau^{(k+2l-1)}$,
and apply Lemma~2.1. 
Then, since ${}^t\!g^{-1}(\tau) = \tau$, we have
$$\split
   & \sum_{r \geq 0} {l \choose r} 
   \big<
   \varXi^{\overline{\overline{k}}}(u) \varTheta^{l-1}
   \omega^r (- 2\rho^*)^{l-r} \tau^m
   \big> \\
    & \quad = \sum_{s \geq 0} \sum_{r \geq 0}
   (-)^{l} {k \choose s} {l \choose r}
   \big<
   \varXi^{\overline{\overline{s}}}(u) 
   \varTheta^{k+l-1-s} 
   \omega^r (2\rho^*)^{l-r} \tau^m
   \big>.
\endsplit
\tag 5.7
$$
Note that 
$\big< 
\varXi^{\overline{\overline{a}}}(u) \varTheta^b \omega^c \rho^{* d} \tau^e
\big>$ is equal to zero, unless $b \ne d$
(compare the order of $e_1,\ldots,e_N$ 
with the order of $e^*_1,\ldots,e^*_N$).
Hence, the left hand side of (5.7) is equal to zero,
unless $r=1$.
Similarly, the right hand side is equal to zero,
unless $s=k$, $r=1$ or $s=k-1$, $r=0$.
Thus we have
$$\align
   &  l 
   \big<
   \varXi^{\overline{\overline{k}}}(u) \varTheta^{l-1}
   \omega (- 2\rho^*)^{l-1} \tau^m
   \big> \\
    & \qquad =  (-)^{l} l
   \big<
   \varXi^{\overline{\overline{k}}}(u) 
   \varTheta^{l-1}
   \omega (2\rho^*)^{l-1} \tau^m
   \big>
   + (-)^{l} k \big< 
   \varXi^{\overline{\overline{k-1}}}(u) 
   \varTheta^{l}
   (2\rho^*)^{l} \tau^m
   \big>.
\endalign
$$
Simplifying this equality, we have
\TagsOnRight
$$
   k \big< \varXi^{\overline{\overline{k-1}}}(u)
   \varTheta^{l} \rho^{* l} \tau^m \big> 
   + l \big< \varXi^{\overline{\overline{k}}}(u)
   \varTheta^{l-1} \rho^{* l-1} \tau^m \omega \big>
   = 0.
\tag"\qed"
$$
\TagsOnLeft

\vskip .10in

\noindent
{\it Proof of Lemma~{\sl 5.12}.}
By the expansion (5.5), we have
$$
   \varXi^{\overline{k}}(u) 
   = \sum_{r \geq 0} (-)^r R^k_r 
   \varXi^{\overline{\overline{k-r}}}(u) \tau^r.
\tag 5.8
$$
Moreover, we have
$$\align
   k W'_{k-1}(u) 
   & = k ( \varXi^{\overline{k-1}}(u) 
   - \frac{k-1}{2} 
   \varXi^{\overline{k-2}}(u) \tau ) \\
   & = \sum_{r \geq 0} (-)^r k R^{k-1}_r 
   \varXi^{\overline{\overline{k-1-r}}}(u) \tau^r 
   - \frac{k(k-1)}{2} \sum_{r \geq 0} (-)^r R^{k-2}_r 
   \varXi^{\overline{\overline{k-2-r}}}(u) \tau^{r+1}.
\endalign
$$
Since $R^k_r = 0$ for $r < 0$, this is equal to
$$\align
   & \sum_{r \geq 0} (-)^r k R^{k-1}_r 
   \varXi^{\overline{\overline{k-1-r}}}(u) \tau^r 
   - \frac{k(k-1)}{2} \sum_{r \geq 0} (-)^{r-1} R^{k-2}_{r-1} 
   \varXi^{\overline{\overline{k-1-r}}}(u) \tau^r \\
   & \qquad = \sum_{r \geq 0} (-)^r 
   (k R^{k-1}_r + \frac{k(k-1)}{2} R^{k-2}_{r-1}) 
   \varXi^{\overline{\overline{k-1-r}}}(u) \tau^r \\
   & \qquad = \sum_{r \geq 0} (-)^r (k-r) R^k_r 
   \varXi^{\overline{\overline{k-1-r}}}(u) \tau^r.
\endalign
$$
Here we used the relation
$$
   k R^{k-1}_r + \frac{k(k-1)}{2} R^{k-2}_{r-1} 
   = (k-2r) R^k_r + r R^k_r
   = (k-r) R^k_r.
$$
Comparing this with (5.8) and applying Lemma~5.11, 
we have the assertion.
\hfill\qed

\vskip .10in

\noindent
{\bf 5.7.}
Combining Lemmas~5.4 and~5.9, we can write $W_k(u)$ as follows.
Using $k R^{k-1}_l = (k-2l) R^k_l$, we have
$$\align
   W_k(u) 
   & = V_k(u) - k V_{k-1}(u) \cdot \tau_{+}
\tag 5.9 \\
\allowdisplaybreak
   & = \sum_{l \geq 0}
       (-)^l R^k_l \varXi^{\overline{k-2l}}(u) \varTheta^l \rho^{* l}
       - k \sum_{l \geq 0} (-)^l R^{k-1}_l 
       \varXi^{\overline{k-2l-1}}(u) \varTheta^l \rho^{* l} \tau_{+} \\
\allowdisplaybreak
   & = \sum_{l \geq 0}
       (-)^l R^k_l \varXi^{\overline{k-2l}}(u) \varTheta^l \rho^{* l}
       - \sum_{l \geq 0} (-)^l (k-2l) R^k_l 
       \varXi^{\overline{k-2l-1}}(u) \varTheta^l \rho^{* l} \tau_{+} \\
\allowdisplaybreak
   & = \sum_{l \geq 0}
       (-)^l R^k_l \varXi^{\overline{k-2l}}(u) \varTheta^l \rho^{* l}
       - \frac{1}{2} \sum_{l \geq 0} (-)^l (k-2l) R^k_l 
       \varXi^{\overline{k-2l-1}}(u) \varTheta^l \rho^{* l} \omega \\
       & \qquad\qquad
       - \frac{1}{2} \sum_{l \geq 0} (-)^l (k-2l) R^k_l 
       \varXi^{\overline{k-2l-1}}(u) \varTheta^l \rho^{* l} \tau \\
\allowdisplaybreak
   & = \sum_{l \geq 0}
       (-)^l R^k_l W'_{k-2l}(u) \varTheta^l \rho^{* l}
       - \frac{1}{2} \sum_{l \geq 0} (-)^l (k-2l) R^k_l 
       \varXi^{\overline{k-2l-1}}(u) \varTheta^l \rho^{* l} \omega \\
\allowdisplaybreak
   & = W'_k(u)
       + \sum_{l \geq 0}
       (-)^{l+1} R^k_{l+1} W'_{k-2l-2}(u) \varTheta^{l+1} \rho^{* l+1} \\
       & \qquad \qquad 
       - \frac{1}{2} \sum_{l \geq 0} (-)^l (k-2l) R^k_l 
       \varXi^{\overline{k-2l-1}}(u) \varTheta^l \rho^{* l} \omega.
\endalign
$$
By Lemma~5.12, we have
$$
   (k-2l-1) \big< W'_{k-2l-2}(u) \varTheta^{l+1} \rho^{* l+1} \big>
   = (l+1) \big< 
   \varXi^{\overline{k-2l-1}}(u) \varTheta^{l} \rho^{*l} \omega 
   \big>.
$$
Multiplying this by 
$ \frac{1}{k-2l-1} R^k_{l+1} = \frac{k-2l}{2l+2} R^k_l$, 
we have
$$
   R^k_{l+1} \big< W'_{k-2l-2}(u) \varTheta^{l+1} \rho^{* l+1} \big>
   = (k-2l) R^k_l 
   \big< 
   \varXi^{\overline{k-2l-1}}(u) \varTheta^{l} \rho^{*l} \omega 
   \big>.
$$
Combining this with (5.9),
we have $\big< W_k(u) \big> =  \big< W'_k(u) \big>$, namely Theorem~5.1.
This finishes the proof of the main theorem.

\vskip .10in

\noindent
{\it Remark.}
This proof of Theorem~4.5 is essentially more difficult 
than that of Theorem~3.8.
In \cite{I4}, we only used simple commutation relations
in $\Lambda_{2N} \otimes U(\frak{o}(S_0))$ to prove Theorem~3.8.
In contrast to this, 
to prove Theorem~4.5, we needed a variable transformation and 
Lemma~2.1 as seen in Section~5.6.

\vskip .20in


%

\enddocument